\documentclass[11pt]{article}

\usepackage{amsmath,makeidx,amssymb,amscd,amsfonts}
\usepackage{epsfig,color,psfrag,booktabs}
\usepackage[left=3cm,top=2.5cm,right=3cm,bottom=3.3cm]{geometry}

\numberwithin{equation}{section}

\newcommand\norm[1]{\|#1\|}
\newcommand\set[1]{\{#1\}}
\newcommand\abs[1]{|#1|}
\newcommand{\f}{\boldsymbol}
\newcommand{\R}{\mathbb R}
\newcommand{\ph}{\varphi}
\newcommand{\xd}{{\f x}}
\newcommand{\yd}{{\f y}}
\newcommand{\td}{{\f t}}
\newcommand{\rd}{{\f r}}
\newcommand{\omd}{{\f \omega}}
\newcommand{\phd}{{\f \varphi}}
\newcommand{\sid}{{\f \sigma}}

\DeclareMathOperator{\A}{\mathcal A}
\DeclareMathOperator{\Mo}{\mathcal M}
\DeclareMathOperator{\Io}{\mathcal I}

\DeclareMathOperator{\Bo}{\mathcal B}

\DeclareMathOperator{\Md}{\mathbf M}
\DeclareMathOperator{\Id}{\mathbf I}
\DeclareMathOperator{\Bd}{\mathbf B}

\newtheorem{theorem}{Theorem}[section]
\newtheorem{lemma}[theorem]{Lemma}
\newtheorem{propo}[theorem]{Proposition}

\newtheorem{remark}[theorem]{Remark}

\def\@abssec#1{\vspace{2ex} {\bfseries #1 } \ignorespaces}
\renewenvironment{abstract}{\@abssec{\abstractname .}}
     {\par\vspace{.1in}}
\newenvironment{keywords}{\@abssec{Keywords:}}
     {\par\vspace{.1in}}
\newenvironment{AMS}{\@abssec{AMS Classifications:}}
     {\par\vspace{.1in}}

\newcommand{\req}[1]{(\ref{eq:#1})}
\allowdisplaybreaks[1]

\begin{document}

\title{On regularization methods of EM-Kaczmarz type}

\author{
M. Haltmeier,%
\thanks{Department of Mathematics, University of Innsbruck,
Technikerstra\ss e 21a, A-6020 Innsbruck, Austria
(\href{mailto:markus.haltmeier@uibk.ac.at}{\tt markus.haltmeier@uibk.ac.at})}
\ \
A.\,Leit\~ao,%
\thanks{Department of Mathematics, Federal University of St. Catarina,
P.O. Box 476, 88040-900 Florian\'opolis, Brazil
(\href{mailto:aleitao@mtm.ufsc.br}{\tt aleitao@mtm.ufsc.br})}
\ \
E.\,Resmerita,%
\thanks{Industrial Mathematics Institute, Johannes Kepler University,
Altenbergerstra\ss e 69, A-4040 Linz, Austria
(\href{mailto:elena.resmerita@ricam.oeaw.ac.at}
{\tt elena.resmerita@ricam.oeaw.ac.at})}
}

\date{\small\today}

\maketitle

\begin{center}
\begin{minipage}{0.95\textwidth}
\abstract{
We consider regularization methods of Kaczmarz type in connection with the
expectation-maximization (EM) algorithm for solving ill-posed equations.
For noisy data, our methods are stabilized extensions of the well established ordered-subsets
expectation-maximization iteration (OS-EM).
We show monotonicity properties of the methods and present a numerical experiment which
indicates that the extended OS-EM methods we propose are much faster than the standard EM
algorithm.}



\end{minipage}
\end{center}

\section{Introduction} \label{sec:intro}

The {\em expectation-maximization} (EM) algorithm  provides approximations
for maximum likelihood estimators of problems with incomplete or noisy data,
which is the usual framework when dealing with inverse or ill-posed problems.
In particular, the EM algorithm for Poisson models is well known for its applications to
astronomical imaging and to PET (positron emission tomography) - see, e.g.
\cite{Ric72}, \cite{VarScheKau85}.

In this work we address inverse problems modeled by operator equations which admit nonnegative solutions,
with the aim of approaching them by combined EM-Kaczmarz strategies.

We begin our study by considering the operator equation
\begin{equation} \label{eq:lin-ip}
\A x  =  y \, ,
\end{equation}
where $\A: L^1(\Omega) \to L^1(\Sigma)$ is a Fredholm integral operator
of the first kind
\begin{equation} \label{eq:def-opA}
(\A x)(s) \ = \ \int_\Omega a(s,t) \, x(t) \, dt, \quad s\in\Sigma \,.
\end{equation}
Nonnegative solutions of \eqref{eq:lin-ip} can be determined by finding minimizers of
the functional
$$
f(x) \ := \ \int_\Sigma \Big[ y(s) \, \ln \frac{y(s)}{(\A x)(s)}
            - y(s) + (\A x)(s) \Big] \, ds \, .
$$
Formally, the first order necessary condition for such a minimizer reads as
\begin{equation} \label{eq:notw-bed}
\A^* \Big( \frac{y}{\A x} \Big) \ = \ \A^* 1 \, .
\end{equation}
If  the assumption $(\A^* 1)(t) = 1$ is satisfied, a solution of \eqref{eq:notw-bed}
can be obtained by solving the corresponding multiplicative fixed-point
equation
\begin{equation} \label{eq:fix-pkt}
    x  \A^* \Big( \frac{y}{\A x} \Big)\ = \ x \,.
\end{equation}
The fixed-point equation \eqref{eq:fix-pkt} motivates the definition of
the EM algorithm, see  \cite{Ric72,VarScheKau85,EggLaR96,Ius92,MueSch87,MueSch89,RedWal84,ResEngIus07},
\begin{equation} \label{eq:def-em}
x_{k+1}(t)
 \ = \
x_k(t)\
\A^* \Big( \frac{y}{\A x_k} \Big)(t)
=
\ x_k(t) \ \int_\Sigma
    \frac{a(s,t) y(s)}{(\A x_k)(s)} \, ds
    \,,
\end{equation}
i.e. an explicit iterative method for solving \eqref{eq:fix-pkt}.
\medskip

The OS-EM (ordered subsets - expectation maximization) iteration was
introduced in \cite{HudLar94} as a computationally more efficient alternative
to the original EM iteration for the discrete case. The main idea is as
follows. The data $y$ are grouped into an ordered sequence of subsets
(or blocks) $y_j$. An iteration of OS-EM consists of a single cycle through
all the subsets, in each subset updating the current estimate by an
application of the EM algorithm in that data subset. This strategy can be
connected to the Kaczmarz type iterative methods recently investigated in
\cite{DecHalLeiSch08, HalLeiSch07, HalKowLeiSch07, KowSch02} for approaching systems of integral
equations.

In order to extend the OS-EM method to infinite dimensional settings, we
first group the data $y$ into $N$ blocks  $y_j := y|_{\Sigma_j}$, where
$\Sigma_j \subset \Sigma$ are not necessarily disjoint and satisfy
$\Sigma = \Sigma_0 \cup \dots \cup \Sigma_{N-1}$. Then equation
\eqref{eq:lin-ip} is decomposed into a system of integral equations of
the first kind
\begin{equation} \label{eq:sys-ip}
    \A_j x \ = \ y_j \,, \qquad j = 0, \dots, N-1 \,,
\end{equation}
where the Fredholm integral operators $\A_j: L^1(\Omega) \to L^1(\Sigma_j)$
correspond to blocks of $\A$ and are defined by
\begin{equation} \label{eq:sys-ops}
(\A_j x)(s) \ := \ \int_\Omega a_j(s,t) \, x(t) \, dt \,,
\end{equation}
with $a_j := a|_{\Omega \times \Sigma_j}$. Notice that $x$ is a solution of
\eqref{eq:sys-ip} if and only if $x$ solves \eqref{eq:lin-ip}.

In order to simplify notation, we drop the indices of the domains
$\Sigma_j$ and simply write $\A_j: L^1(\Omega) \to L^1(\Sigma)$
and $y_j \in L^1(\Sigma)$. Thus, the system of integral equations
\eqref{eq:sys-ip} can be approached by simultaneously minimizing
$$
f_j(x) \ := \ \int_\Sigma \Big[ y_j(s) \, \ln \frac{y_j(s)}{(\A_j x)(s)}
              - y_j(s) + (\A_j x)(s) \Big] \, ds \, .
$$
It is worth noticing that $f_j(x) = d(y_j, \A_j x)$, where $d(u,v)$ is the
Kullback-Leibler (KL) distance defined by
\begin{equation} \label{eq:def-KLd}
d(v,u) \ := \ \int \Bigl[ v(t)  \ \ln \frac{v(t)}{u(t)} - v(t) + u(t)\Bigr] \, dt \, .
\end{equation}
Throughout this article we will make use of the KL-distance $d(v,u)$ with either $u, v \in L^1(\Omega)$ or
$u, v \in  L^1(\Sigma)$.

\begin{remark} \label{rem:adj-cond}
Analog as in \eqref{eq:notw-bed}, if the assumption $\A_j^* \, 1 = 1$, $j = 0, \dots, N-1$, is satisfied, then
the first order necessary condition for
a minimizer of $f_j$ is given by $\A_j^* \big( y_j / (\A_j x) \big) = 1$,
and the corresponding multiplicative fixed-point equation reads
$P_j (x) := x \A_j^* \big( y_j / (\A_j x) \big) = x $.
\end{remark}

The OS-EM algorithm corresponds to a Kaczmarz type method for solving
system \eqref{eq:sys-ip} and can be written in the form
\begin{equation} \label{eq:def-osem}
x_{k+1}
=
P_j (x_k)
=
x_k  \int_\Sigma
           \frac{a_j(s,\cdot) \, y_j(s)}{(\A_j x_k)(s)} \, ds
\, ,
\end{equation}
where the index $0 \le j < N$ relates to the iteration index $k$ by the
formula $j = [k] := (k$ mod $N)$.  Clearly, the case $N=1$ corresponds to
the standard EM algorithm.

The cyclic structure of the iteration in \eqref{eq:def-osem} is easily
recognizable (each cycle consists of $N$ steps).
Notice that each step within a cycle is an explicit step for solving the
fixed point equation $x \A_{[k]}^* \big( y_{[k]}
/(\A_{[k]} x) \big) = x$,
and can be interpreted as an EM iterative step for solving the $[k]$-th
equation (or block) of system \eqref{eq:sys-ip}.
\bigskip

This article is outlined as follows.
In Section~\ref{sec:2} we formulate a series of assumptions, which are
necessary for the analytical investigation of the OS-EM method.
Moreover, we present some basic results concerning the KL-distance.
Section~\ref{sec:3} contains an analysis of the OS-EM iteration
\eqref{eq:def-osem}, i.e., monotonicity results and consequences concerning
the asymptotic behavior of the iterations.
Section~\ref{sec:4} studies the case of noisy data and introduces the
{\em loping OS-EM method} \eqref{eq:def-osem-noise} which is a modification
of the OS-EM iteration for noisy data.
Stability results that use discrepancy type
principles are stated.
In Section~\ref{sec:5} we present some numerical experiments regarding
application of the OS-EM methods to the inversion of the circular
Radon transform.
Section~\ref{sec:6} is devoted to final remarks and conclusions.

\section{Assumptions and basic results} \label{sec:2}

Throughout this article we assume the domains $\Omega$ and $\Sigma$
in Section~\ref{sec:intro} to be open bounded subsets of ${\mathbb R}^d$, $d \ge 1$. The parameter space for investigating system \eqref{eq:sys-ip} is
\begin{equation} \label{eq:def-delta}
\Delta \ := \ \{ x \in L^1(\Omega) \, ; \ x \ge 0 \, , \
                 \int_\Omega x(t) \, dt \, = \, 1 \} \, ,
\end{equation}
and the starting element $x_0$ of iteration \eqref{eq:def-osem} is chosen
such that $x_0 \in \Delta$.

Moreover, we make the following assumptions to the framework introduced in
Section~\ref{sec:intro}:
\begin{itemize}
\item[(A1)]
The kernel functions $a_j: \Sigma \times \Omega \to \mathbb R$, $j = 0,\dots,N-1$, in
\eqref{eq:def-opA} satisfy $\int_\Sigma a_j(s,t) \, ds = 1$ for a.e.
$t \in \Omega$;
\item[(A2)]
There exist positive constants $m$ and $M$ such that
$m \le a_j(s,t) \le M$ a.e. in $\Sigma \times \Omega$;
\item[(A3)]
The exact data $y_j \in L^1(\Sigma)$ in \eqref{eq:sys-ip}
satisfy $\int_\Sigma y_j(s) \, ds = 1$; moreover, there exists $M' > 0$
such that $y_j(s) \le M'$ a.e. in $\Sigma$;%
\item[(A4)]
System \eqref{eq:sys-ops} has a non-negative solution
$x^* \in L^1(\Omega)$, which does not vanish a.e. in $\Omega$; moreover,
$d(x^*, x_0) < \infty$.
\end{itemize}

Assumption  (A2) implies that the operators $\A_j: L^1(\Omega) \to L^1(\Sigma)$ are continuous.
Moreover, any $\A_j x_k$ is in $L^{\infty}(\Sigma)$ and
bounded away from zero. This further ensures that $1/\A_j x_k$ has the same
properties and then yields that the integrals in \eqref{eq:def-osem} are
well-defined.

\medskip
In the sequel we discuss some basic properties of the KL-distance in
\eqref{eq:def-KLd} that will be needed in the forthcoming sections.
This functional plays a key role in the convergence analysis of the
OS-EM method. For details, we refer the reader to \cite{ResEngIus07, ResAnd07}.

\begin{lemma} \label{lem:kl-prop}
Let $u$ and $v$ be two $L^1$ functions such that $(u,v)$ is in the domain  of the
KL-distance $d(\cdot,\cdot)$ defined in \eqref{eq:def-KLd}.
The following assertions hold true:
\begin{itemize}
\item[(i)]   $d(v,u) \ge 0$ and $d(v,u) = 0$ iff $v = u$ a.e.;
\item[(ii)]  $\| v - u \|_{L^1}^2 \le \big( \frac{2}{3} \| v \|_{L^1} +
             \frac{4}{3} \| u \|_{L^1} \big) \, d(v,u)$;
\item[(iii)] The function $(v,u) \mapsto d(v,u)$ is convex;
\item[(iv)]  Let $\{ v_n \}$ and $\{ u_n \}$ be given sequences in $L^1$.
If $\{ u_n \}$ is bounded and  $\lim\limits_{n\to\infty} d(v_n,u_n) = 0$,
then $\lim\limits_{n\to\infty} \| v_n - u_n \|_{L^1} = 0$.
\end{itemize}
\end{lemma}

\section{The OS-EM method  for exact data} \label{sec:3}

The first result of this section relates to a monotonicity property of the
OS-EM iteration.

\begin{lemma} \label{lem:monot1}
Let Assumptions (A1)-(A3) be satisfied, let $x \in \Delta$, and
denote $P_j (x) = x \A_j^*(y_j/\A_ j x )$.
Then the following assertions hold true:
\begin{itemize}
\item[(i)]  $P_j(x) \in \Delta$ and
            $d(P_j(x), x) \le f_j(x) - f_j(P_j(x))$,
            for $j = 0, \dots, N-1$;

\item[(ii)] If $x^* \in \Delta$ is a minimizer of $f_j$ for some $0 \le j \le
            N-1$, and $d(x^*, x) < \infty$, then $d(x^*, P_j(x)) < \infty$ and
            \ $f_j(x) - f_j(x^*) \le d(x^*, x) - d(x^*, P_j(x))$.
\end{itemize}
\end{lemma}
\begin{proof}
Results immediately from \cite[Prop.~3.1]{ResEngIus07} applied to the function
$f_j$ and the corresponding $P_j$. \hfill
\end{proof}

\medskip
From Lemma~\ref{lem:monot1}~(i) and Lemma~\ref{lem:kl-prop}~(i) we conclude
that $f_j(P_j(x)) \le f_j(x)$ and $P_j(x) \in \Delta$.
Moreover, if $x^* \in \Delta$ is a solution of \eqref{eq:sys-ip} with
$d(x^*, x) < \infty$, then   $x^*$ minimizes $f_j$, for every $j = 0, \dots, N-1$.
Lemma~\ref{lem:monot1}~(ii) and the fact that $f_j(x^*)=0$ therefore yield
\begin{equation} \label{eq:rem-monot}
f_j(x) \ \le \ d(x^*, x) - d(x^*, P_j(x)) \, , \quad
j = 0, \dots, N-1 \, .
\end{equation}

In the next lemma we reinterpret the inequalities derived in Lemma~%
\ref{lem:monot1} in terms of the OS-EM iteration.

\begin{lemma} \label{lem:monot2}
Let Assumptions (A1)-(A3) be satisfied, and let $\{ x_k \}$ be defined by
iteration \eqref{eq:def-osem}. Then the following assertions hold true:
\begin{itemize}
\item[(i)]  $d(x_{k+1}, x_k) \, \le \, f_{[k]}(x_k) - f_{[k]}(x_{k+1})$;
\item[(ii)] If $x^* \in \Delta$ is a solution of \eqref{eq:sys-ip}, then
$f_{[k]}(x_{k}) \, \le \, d(x^*, x_k) - d(x^*, x_{k+1})$.
\end{itemize}
\end{lemma}
\begin{proof}
Results from Lemma \ref{lem:monot1} and \eqref{eq:rem-monot}. \hfill
\end{proof}
\medskip

In the next theorem we formulate the main monotonicity results for the
OS-EM iteration with respect to the KL-distance, as well as convergence
results in case the iterations are bounded.

\begin{theorem} \label{th:monot}
Let Assumptions (A1)-(A3) be satisfied, and the sequence $\{ x_k \}$ be
defined by iteration \eqref{eq:def-osem}.
Then we have
\begin{itemize}
\item[(i)]  $f_{[k]}(x_{k+1}) \le f_{[k]}(x_k)$, for every $k \in \mathbb N$.
\end{itemize}

\noindent
Moreover, if assumption (A4) is satisfied, then the following assertions hold
true:

\begin{itemize}
\item[(ii)]  The sequence $\{ d(x^*, x_k) \}$ is nonincreasing;
\item[(iii)]  $\lim\limits_{k\to\infty} f_{[k]}(x_k) = 0$;
\item[(iv)]   $\lim\limits_{k\to\infty} d(x_{k+1}, x_k) = 0$;
\item[(v)]  For each $0 \le j \le N-1$ and $p\in[1,\infty)$ we have
             \begin{equation} \label{eq:l1-conv}
             \lim_{m\to\infty} \| \A_j x_{j+mN} - y_j \|_{L^p(\Sigma)} = 0 \, .
             \end{equation}

\item[(vi)] If  $\{x_k\}$ is bounded in some
$L^p(\Omega)$ space, with $p\in (1, \infty)$, then it has a subsequence
which converges weakly in $L^p(\Omega)$ to  a solution of system  \eqref{eq:sys-ip}.
\end{itemize}
\end{theorem}
\begin{proof}
Items (i) and (ii) follow from Lemma~\ref{lem:monot2}~(i) and Lemma~%
\ref{lem:kl-prop}~(i).
Item (ii) implies the existence of $\mu \ge 0$
such that $\lim\limits_{k\to\infty} d(x^*, x_k) = \mu$. Thus, (iii)
follows from Lemma~\ref{lem:monot2}~(ii).

To prove (iv), notice that (i) and~(iii) imply
\begin{equation} \label{eq:monot-aux}
\lim_{k\to\infty} f_{[k]}(x_{k+1}) = 0 \, .
\end{equation}
Now, (iv) results from \eqref{eq:monot-aux},  Item (iii) and Lemma~%
\ref{lem:monot2}~(i).

Next we prove (v).
Since $f_j(x) = d(y_j, \A_j x)$, it follows from
(iv) that $\lim\limits_{k\to\infty} d( y_{[k]} , A_{[k]} x_k) = 0$. Consequently,
$\lim\limits_{m\to\infty} d(y_j , \A_j x_{j+mN}) = 0$,  for every  $j = 0, \dots, N-1$.
Now,  by applying Lemma~\ref{lem:kl-prop}~(iv) we obtain \eqref{eq:l1-conv} for $p=1$.
The case $p \in (1, \infty)$ follows
from \cite[Prop.~4.1 and Lem.~4.2]{ResEngIus07}.

\smallskip \noindent
The proof of assertion (vi) is divided in several parts:

\smallskip \noindent
(1)
Claim:
$x_k \in L^\infty(\Omega)$ for each $k \in \mathbb N$.

\noindent
Since $x_0 \in \Delta$ by hypothesis, we have $m \le (A_0 x_0)(s) \le M$
a.e. in $\Sigma$ by (A2). Consequently, it results from (A2) and (A3)
that \[\frac{a_1(s,t) y_0(s)}{(A_0 x_0)(s)} \le \frac{M M'}{m} \,,  \qquad \text{ a.e. in }
\Sigma \times \Omega\,,
\]
and from \eqref{eq:def-osem} follows $x_1 \in L^\infty(\Omega)$.
Part (1) follows by induction if one
observes that $x_k \in L^\infty(\Omega)$ together with (A2) and (A3)
imply $x_{k+1} \in L^\infty(\Omega)$.

\smallskip \noindent
(2)
By hypothesis, the sequence $\{x_k\}$ is bounded
in $L^p(\Omega)$.
Therefore, there is a  subsequence  denoted again by
$\{ x_k \}_{k\in\mathbb N}$, which converges weakly  in $L^p(\Omega)$
to some $z \in L^p(\Omega)$, for some $p\in(1,+\infty)$.

\smallskip \noindent
(3)
Conclusion of the proof under a simplifying assumption.

\noindent
Let us assume for the moment that, for each fixed $0 \le j \le N-1$,
the subsequence $\{ x_k \}_{k\in\mathbb N}$ obtained in Part~(2)
contains infinitely many indices of the form $k = j + m  N$, $m  \in \mathbb N$.
Then, for $j=0$, we can extract from $\{ x_k \}$ a subsequence
$\{ x_{k_i} \}$ with indices of the form $k_i = m_i N$.
Obviously, $[k_i] = 0$ for all indices of the subsequence $\{ x_{k_i} \}$,
and from (vi) it follows that $\A_0 x_{k_i} \to y_0$ strongly, and
thus weakly in $L^p(\Sigma)$.
Since $\A_0$ is continuous from $L^p(\Omega)$ to $L^p(\Sigma)$ due
to (A2), it is weakly continuous. Part~(2) implies that
$\A_0 x_{k_j} \to \A_0 z$ weakly in $L^p(\Sigma)$. From the uniqueness of
weak limits it follows that $\A_0 z = y_0$.

\noindent
By repeating the argumentation for
$j = 1, \dots, N-1$, we conclude that $\A_j z = y_j$, for $j = 1, \dots, N-1$,
thus proving that $z$ is a solution of system \eqref{eq:sys-ip}.
\smallskip

\smallskip \noindent
(4)
Conclusion of the proof in the general case.

\noindent
If the assumption in part (3) does not hold, then there must be at least
one $0 \le j_0 \le N-1$ such that the subsequence $\{ x_k \}_{k\in\mathbb N}$
obtained in part (3) contains infinitely many indices of the form
$k = j_0 + mN$, $m \in \mathbb N$.
Arguing as in part~(4) we obtain a subsequence $\{ x_{k_i} \}$, with indices
of the form $k_i = j_0 + m_i N$, $i \in \mathbb N$, such that $A_{j_0} x_{k_i} \to
y_{j_0}$ weakly in $L^p(\Sigma)$, and conclude that $A_{j_0} z = y_{j_0}$.

\noindent
Now, let us consider the subsequence $\{ x_{k_i+1} \}$, with indices of the
form $k_i+1 = (j_0+1)+m_iN$, $i \in \mathbb N$.%
\footnote{Notice that $\{ x_{k_i+1} \}$ may contain elements which do not
belong to the convergent subsequence $\{ x_k \}$ obtained in part~(3),
while $\{ x_{k_i} \}$ is a subsequence extracted from $\{ x_k \}$.}
Item (iv) implies that $d(x_{k_i+1}, \, x_{k_i}) \to 0$, as
$i\to\infty$.
Since both subsequences $\{ x_{k_i+1} \}_i$, $\{ x_{k_i} \}_i$ are in $L^p(\Omega)$
and bounded (part (1) above), it follows from Lemma~\ref{lem:kl-prop}~(iv)
that $\| x_{k_i+1} - x_{k_i} \|_{L^p} \to 0$ as $k_i \to \infty$.
Therefore, $x_{k_i+1} \to z$ weakly in $L^p(\Omega)$%
\footnote{Notice that $x_{k_i} \to z$ weakly in $L^p(\Omega)$.}
and from the continuity of $\A_{[j_0+1]}$ follows that $\A_{[j_0+1]} x_{k_i+1}
\to \A_{[j_0+1]}z$ weakly in $L^p(\Sigma)$.
Moreover,  (v) yields $\A_{[j_0+1]} x_{k_i+1}
\to y_{j_0+1}$ weakly in $L^p(\Sigma)$. Then we conclude that $\A_{[j_0+1]} z
= y_{[j_0+1]}$.

\noindent
Repeating the argumentation for the subsequences $\{ x_{k_i+2} \}$, \dots,
$\{ x_{k_i+N-1} \}$, we conclude that $\A_j z = y_j$, for every $j$,
proving that $z$ is a solution of system \eqref{eq:sys-ip}.
\end{proof}
\bigskip

\begin{remark} \label{rem:conv-sseq}
We can interpret Theorem~\ref{th:monot}~(v) as follows:
If we consider the  subsequence $\{ x_{j+mN} \}_{m\in\mathbb N}$ formed by the
$j$-th component of each cycle of the OS-EM iteration (where
$0 \le j \le N-1$), then the $L^1$-norm of the residual corresponding to
this subsequence converges to zero.

Moreover, Theorem~\ref{th:monot} (iv) guarantees that, given any
two "consecutive" subsequences $\{x_{j + mN}\}_{m\in\mathbb N}$ and
$\{x_{(j+1)+mN}\}_{m\in\mathbb N}$, we have
$$
\lim\limits_{m\to\infty} d(x_{(j+1)+mN}, \, x_{j+mN}) \ = \ 0 \, ,
$$
for each $j = 0, \dots, N-2$.
\end{remark}

\section{The loping OS-EM method for noisy data} \label{sec:4}

Our next goal is to modify the OS-EM iteration by introducing a
relaxation parameter, and to investigate monotonicity and stability results
for this modified iteration (the so called {\em loping OS-EM method})
in the case of noisy data. As remarked in \cite{HudLar94},
``With noisy data though, inconsistent applications (of discrete OS-EM
-- authors' note) result.''

We aim at characterizing the loping OS-EM method as an
{\em iterative regularization method} in the sense of \cite{EngHanNeu96}.

For the rest of this section we assume that the right hand side
of \eqref{eq:sys-ip} is not exactly known. Instead, we have only
approximate measured data $y_j^\delta \in L^1(\Sigma)$ satisfying
\begin{equation} \label{eq:noisy-data}
\|y_j-y_j^{\delta}\|_{L^1} \ \le \ \delta_j \, , \qquad j = 0, \dots, N-1 \, .
\end{equation}
We denote $\delta := (\delta_0, ..., \delta_{N-1})$.

In this noisy data case we are interested in finding an approximate
solution for the system
\begin{equation} \label{eq:sys-ip-noise}
\A_j x \ = \ y_j^\delta \, , \qquad j = 0, \dots, N-1 \,.
\end{equation}
The following assumptions are required for the analysis:
\begin{itemize}
\item[(A5)] The noisy data $y_j^\delta \in L^1(\Sigma)$ satisfies
$\int_\Sigma y_j^\delta(s) \, ds = 1$.

\item[(A6)] There exist $M_1, m_1 > 0$ such that $M_1 \geq y_j^\delta \ge m_1$
a.e. in $\Sigma$.
\end{itemize}

Also necessary for the analysis are the following functions
associated to the equations of system \eqref{eq:sys-ip-noise}
\begin{equation}\label{eq:fdelta}
f_j^\delta(x) \ := \ \int_\Sigma \Big[ y^\delta_j(s) \,
                  \ln \frac{y_j^\delta(s)}{(\A_j x)(s)}
                  - y_j^\delta(s) + (\A_j x)(s) \Big] \, ds
                  \,.
\end{equation}
Notice that $f_j^\delta(x) = d(y_j^\delta, \A_j x)$.

The {\em loping OS-EM iteration} for the inverse problem
\eqref{eq:sys-ip-noise} with noisy data is defined by
\begin{subequations} \label{eq:def-osem-noise}
\begin{equation} \label{eq:def-osem-noise1}
x_{k+1}^\delta \ = \ x_k^\delta \, \omega_k
\end{equation}
where
\begin{equation} \label{eq:def-osem-noise2}
\omega_k \ = \
  \begin{cases}
    \int_\Sigma \frac{a_{[k]}(s,\cdot) \, y_{[k]}^\delta(\cdot)}
    {(A_{[k]} x_k^\delta)(s)} \, ds \ =: \ P_{[k]}^\delta (x_k^\delta)\,,
    & \ f_{[k]}^\delta(x_k^\delta) > \tau \gamma \delta_{[k]} \\
    1\,, & \ \mbox{ else}
  \end{cases} .
\end{equation}
\end{subequations}
The constants $\tau$ and $\gamma$ in \eqref{eq:def-osem-noise2} are chosen
such that
\begin{equation} \label{eq:def-tau-gamma}
\tau \ > \ 1 \, , \quad\quad
\gamma \ = \max \Big\{ \big| \ln \frac{m_1}{M} \big| \, , \,
                       \big| \ln \frac{M_1}{m} \big| \Big\} \, ,
\end{equation}
where $m$, $M$, $m_1$, $M_1$ are the positive constants defined in (A2) and (A6).

\begin{remark}
It is worth noticing that, for noisy data, the iteration in
\eqref{eq:def-osem-noise} is much different from the iteration in
\eqref{eq:def-osem}: The relaxation parameter $\omega_k$ effects
that the iterates defined in \eqref{eq:def-osem-noise1} become
stationary if all components of the residual vector
$d(y_{[k]}^\delta, A_{[k]} x_k^\delta)$ fall below a pre-specified threshold.

Another consequence of using these relaxation parameters is the fact that,
after a large number of iterations, $\omega_k = 1$ for some $k$ within each
iteration cycle. Therefore, the computational evaluation of the adjoint
operator
$$
\A_{[k]}^* \Big( \frac{y_{[k]}^\delta}{A_{[k]} x_{k}^\delta} \Big) \ = \
\int_\Sigma \frac{a_{[k]}(s,\cdot) \, y_{[k]}^\delta(\cdot)}
                {(\A_{[k]} x_k^\delta)(s)} \, ds
$$
might be loped, making the loping OS-EM iteration in
\eqref{eq:def-osem-noise} a fast alternative to the OS-EM method.

In the case of noise free data, i.e. $\delta_j = 0$ in \eqref{eq:noisy-data},
we choose $\omega_k = P_{[k]}^\delta(x_k^\delta) = P_{[k]}(x_k)$ and the
loping OS-EM iteration \eqref{eq:def-osem-noise} reduces to the OS-EM
method \eqref{eq:def-osem}.
\end{remark}

In the sequel we prove a monotonicity result for the loping OS-EM iteration
in the case of noisy data. First however, we derive an auxiliary estimate.

\begin{lemma} \label{lem:monot-noise}
Let assumptions (A1)-(A5) hold true. Moreover, let $y_j^\delta$, $\delta_j$
be given as in \eqref{eq:noisy-data}, with $\delta_{j_0} > 0$ for some
$0 \le j_0\le N-1$. Then we have
\begin{equation} \label{eq:monot-aux-noise}
f_{[k]}^\delta(x_k^\delta) - d(y_{[k]}, y_{[k]}^\delta) \ \le \
d(x^*, x_k^\delta) - d(x^*, x_{k+1}^\delta) \, .
\end{equation}
for all $k \in \mathbb N$ with $[k] = j_0$.
\end{lemma}
\begin{proof}
Since (A1)-(A3) are satisfied, we argue as in the proof of
\cite[Prop.~5.2]{ResEngIus07} to conclude that for every $v$, $w \in \Delta$,
and $0 \le j \le N-1$ the inequality
$$
d(P_j(w), P_j^\delta(v)) \le d(y_j, y_j^\delta) + d(P_j(w), v) - d(P_j(w), w)
                             + f_j(w) - f_j(v) \, ,
$$
holds true. Therefore, given $k \in \mathbb N$ with $[k] = j_0$,
\eqref{eq:monot-aux-noise} follows by taking $j = [k]$, $w = x^*$,
$v = x_k^\delta$, and by observing that $P_{[k]}(x^*) = x^*$. \hfill
\end{proof}

\begin{propo} \label{pr:monot-noise}
Let assumptions (A1)-(A6) hold true and $\tau$, $\gamma$ be defined as
in \eqref{eq:def-tau-gamma}. Moreover, let $y_j^\delta$, $\delta_j$ be given
as in \eqref{eq:noisy-data} with $\delta_j > 0$ for $j = 0, \dots, N-1$.
Then the sequence $\{ x_k^\delta \}$ defined by iteration
\eqref{eq:def-osem-noise} satisfies
\begin{equation} \label{eq:monot-noise}
d(x^*, x_{k+1}^\delta) \ \le \ d(x^*, x_k^\delta) \,,
\qquad
\ k \in {\mathbb N} \, .
\end{equation}
\end{propo}
\begin{proof}
If $f_{[k]}^\delta(x_k^\delta) \le \tau \gamma \delta_{[k]}$, then $w_k = 1$
by \eqref{eq:def-osem-noise2}. Therefore, $x_{k+1}^\delta = x_k^\delta$ and
\eqref{eq:monot-noise} follows with equality.
If $f_{[k]}^\delta(x_k^\delta) > \tau \gamma \delta_{[k]}$, notice that a simple
calculation yields
$$
d(x^*, x_k^\delta) - d(x^*, x_{k+1}^\delta)
\ \ge \
f_{[k]}^\delta(x_k^\delta) + \int_\Sigma [y_{[k]}(s) - y_{[k]}^\delta(s)] \, \ln
        \Big( \frac{y_{[k]}^\delta(s)}{(A_{[k]} x_k^\delta)(s)} \Big) \, ds
$$
from \eqref{eq:monot-aux-noise}.
Therefore, \eqref{eq:monot-noise} follows from
\begin{align}
f_{[k]}^\delta(x_k^\delta) + \int_\Sigma [y_{[k]} - y_{[k]}^\delta]
& \ln
\Big( \frac{y_{[k]}^\delta}{A_{[k]} x_k^\delta} \Big) \, ds \
\ge
\nonumber \\
& \ge f_{[k]}^\delta(x_k^\delta) - \| y_{[k]} - y_{[k]}^\delta \|_{L^1} \,
        \| \ln \Big( \frac{y_{[k]}^\delta}{A_{[k]} x_k^\delta} \Big) \|_{L^\infty}
        \nonumber \\
& \ge f_{[k]}^\delta(x_k^\delta) - \delta_{[k]}
        \max \Big\{ \big| \ln \frac{m_1}{M} \big| \, , \,
                    \big| \ln \frac{M_1}{m} \big| \Big\} \nonumber \\[1ex]
& \ge f_{[k]}^\delta(x_k^\delta) - \gamma \, \delta_{[k]}
        \label{eq:monot-noise-aux} \\[1.5ex]
& \ge (\tau - 1) \, \gamma \, \delta_{[k]} \nonumber
\end{align}
together with \eqref{eq:def-tau-gamma}. To obtain the inequalities above we
used \eqref{eq:noisy-data}, \eqref{eq:def-tau-gamma}, (A2) and (A6). \hfill
\end{proof}
\medskip

Proposition~\ref{pr:monot-noise} gives us a hint on how to choose the stopping
rule for the loping OS-EM iteration. That is, we stop the iteration at
\begin{equation} \label{eq:def-discrep1}
k_*^\delta  :=  {\min} \{ mN \in {\mathbb N}\, ; \ x_{mN}^\delta =
               x_{mN+1}^\delta = \cdots = x_{mN+N-1}^\delta \} \, .
\end{equation}
In other words, $k_*^\delta$ is the smallest integer multiple of $N$ such that
\begin{equation} \label{eq:def-discrep2}
x_{k_*^\delta} = x_{k_*^\delta+1} = \dots = x_{k_*^\delta+N-1} \, .
\end{equation}

In the sequel, we prove that the stopping index $k_*^\delta$ in
\eqref{eq:def-discrep1} is well defined and that the corresponding iterations
stably converge to a solution of the system, if they are bounded in some
$L^p$ space with $p\in (1,+\infty)$.

\begin{theorem} \label{pro:stop-rule}
Let assumptions~(A1)-(A6) be satisfied, and $k_*^\delta \in \mathbb N$ be
chosen according to \eqref{eq:def-discrep1}. Then the following assertions
hold true:
\begin{itemize}
\item[(i)]   The stopping index $k_*^\delta$ defined in \eqref{eq:def-discrep1}
             is finite;
\item[(ii)]  More precisely, $k_*^\delta = O(\mathbf\delta_{\min}^{-1})$, where
             $\mathbf\delta_{\min} := \min \{ \delta_0, \dots, \delta_{N-1}\}$;

\item[(iii)]
$d(y_j^\delta, \A_j x_{k_*^\delta}^\delta)  \leq  \tau \gamma \delta_j$,
for every $j = 0, \dots, N-1$;

\item[(iv)]
For every $p\in[1,+\infty)$ and every $j = 0, \dots, N-1$ we have
             \[
             \lim_{\delta \to 0} \norm{\A_{j} x_{k_*^\delta}^\delta - y_j}_{L^p(\Sigma)} = 0 \,.
             \]

\item[(v)] Let
$\{ \delta^l := (\delta^l_{0}, \dots, \delta^l_{N-1}) \}_{l\in\mathbb N}$ be
a sequence in $(0,\infty)^N$ with $\lim\limits_{l\to\infty} \delta^l_j = 0$,
for each $0\le j \le N-1$.
Moreover, let $\{ y^l := (y^l_{0}, \dots, y^l_{N-1}) \}_{l\in\mathbb N}$ be a
sequence of noisy data satisfying
\begin{equation} \label{eq:assump-stabil-th}
\|y_j-y_j^l\|_{L^1} \ \le \ \delta^l_j \, , \ \ j = 0, \dots, N-1 \, ,
\ \ l \in \mathbb N \, ,
\end{equation}
and $k_*^l := k_*^{\delta^l} = k_*(\delta^l, y^l)$ denote the corresponding
stopping index defined in \eqref{eq:def-discrep1}. If the sequence
$\{ x_{k_*^l}^{\delta^l} \}_{l\in\mathbb N}$ is bounded in some $L^p(\Omega)$
space, with $p\in (1,+\infty)$, then it has a subsequence which converges
weakly in $L^p(\Omega)$ to a solution of system  \eqref{eq:sys-ip}.
\end{itemize}
\end{theorem}
\begin{proof}
(i) Assume by contradiction that $k_*^\delta$ is not finite. Then it
results from \eqref{eq:def-discrep1} that $x_{k+1}^\delta \not=
x_k^\delta$ at least once in each cycle of iteration
\eqref{eq:def-osem-noise}. Hence for every $m \in \mathbb N$
there exits $j_m \in \{ 0, \dots, N-1 \}$ such that
\begin{equation} \label{eq:argument-j0}
f_{j_m}^\delta(x_{j_m+mN}^\delta) > \tau \gamma \delta_{j_m} \,.
\end{equation}
From \eqref{eq:monot-noise-aux} in the proof of Proposition~\ref{pr:monot-noise},
it follows
$$
d(x^*, x_k^\delta) - d(x^*, x_{k+1}^\delta) \ \ge \
\max\{ f_{[k]}^\delta(x_k^\delta) - \gamma \delta_{[k]} , \, 0 \} \, , \quad
k \in \mathbb N \, .
$$
Summing up this inequality for $k = 0, \dots, lN-1$ implies
\footnote{Notice that $x_0^\delta = x_0$.}
\begin{multline*}
d(x^*, x_0) - d(x^*, x_{lN}^\delta)
\geq
\sum_{k=0}^{lN-1}
\max \{ f_{[k]}^\delta(x_k^\delta) - \gamma\delta_{[k]}
, \, 0 \} \\
= \sum_{m=0}^l \sum_{j=0}^{N-1} \max\{
f_{j}^\delta(x_{j+mN}^\delta) - \gamma\delta_{j} , \, 0 \} \, ,
\quad l \in \mathbb N \,.
\end{multline*}
Then, it follows from \eqref{eq:argument-j0}
\begin{equation} \label{eq:k*-finite}
d(x^*, x_0) \ \ge \
\sum_{m=0}^l \left( f_{j_m}^\delta(x_{j_m+mN}^\delta) - \gamma\delta_{j_m} \right) \\[-1ex]
> \sum_{m=0}^l (\tau - 1) \, \gamma \, \delta_{j_m}
> l \ (\tau - 1) \, \gamma \, \delta_{\min}
\,, \quad l \in \mathbb N \,.
\end{equation}
However, due to \eqref{eq:def-tau-gamma}, the right hand side of
\eqref{eq:k*-finite} becomes unbounded as $l\to\infty$, contradicting (A4).
Therefore, $k_*^\delta$ must be finite.
To prove (ii), it is enough to take $l = k_*^\delta / N \in \mathbb N$ in
\eqref{eq:k*-finite} and obtain $k_*^\delta  < N d(x^*, x_0)/\bigl( (\tau-1) \gamma\delta_{\rm  min}\bigr)$.

To prove (iii), we assume by contradiction that
\[ f_{j_0}^\delta(x_{k_*^\delta}^\delta) =
d(y_{j_0}^\delta,  \A_{j_0} x_{k_*^\delta}^\delta) > \tau \gamma \delta_{j_0}\,,\]
for some $0 \le j_0 \le N-1$. Thus, it results from \eqref{eq:def-discrep2}
that $f_{j_0}^\delta(x_{k_*^\delta+j_0}^\delta) > \tau \gamma \delta_{j_0}$.
Therefore, it follows from \eqref{eq:monot-noise-aux} in the proof of
Proposition~\ref{pr:monot-noise} that
$$
0 \  =  \ d(x^*, x_{k_*^\delta+j_0}^\delta) - d(x^*, x_{k_*^\delta+j_0+1}^\delta)
  \ \ge \ f_{j_0}^\delta(x_{k_*^\delta+j_0}^\delta) - \gamma \delta_{j_0}
  \ \ge \ (\tau - 1) \gamma \delta_{j_0} \, ,
$$
which contradicts \eqref{eq:def-tau-gamma}.

(iv) and (v) The proofs follow the lines  of the proof of Theorem \ref{th:monot}~(v), (vi). \hfill
\end{proof}

\begin{remark}[Stability for noisy data in $L^2(\Sigma)$]
When dealing with inverse problems, bounds for the noisy data are
most commonly given in the $L^2$-norm, i.e. the approximate measured
data $y_j^\delta \in L^2(\Sigma)$ is assumed to satisfy
\begin{equation} \label{eq:noisy-data2}
\| y_j - y_j^\delta \|_{L^2} \ \le \ \delta_j \, , \quad j = 0, \dots, N-1 \, ,
\end{equation}
instead of \eqref{eq:noisy-data}. In this case, the loping OS-EM iteration
is defined by \eqref{eq:def-osem-noise}, where the ``loping condition"
$f_{[k]}^\delta(x_k^\delta) > \tau \gamma \delta_{[k]}$ in
\eqref{eq:def-osem-noise2} is substituted by
\begin{equation} \label{eq:lop-cond-l2}
f_{[k]}^\delta(x_k^\delta) \ > \
\tau \delta_{[k]} \| \ln(y_{[k]}^\delta / (\A_{[k]} x_k^\delta)) \|_{L^2} \, .
\end{equation}
Under this assumptions it is possible to state a stability result,
similar to the one  in Theorem \ref{pro:stop-rule} (iv). One
argues as follows:
\begin{itemize}
\item First of all, notice that monotonicity of the error with respect to
the KL-distance (as in \eqref{eq:monot-noise}) follows when using the
Cauchy-Schwarz inequality in $L^2(\Sigma)$ to derive the estimate (compare
with \eqref{eq:monot-noise-aux})
\begin{equation*}
f_{[k]}^\delta(x_k^\delta) + \int_\Sigma [y_{[k]} - y_{[k]}^\delta] \,
\ln (y_{[k]}^\delta / (\A_{[k]} x_k^\delta)) \, ds
\ge
(\tau - 1) \, \delta_{[k]} \| \ln(y_{[k]}^\delta / (\A_{[k]} x_k^\delta)) \|_{L^2}
\,.
\end{equation*}
\item By defining the stopping index $k_*^\delta$ as in \eqref{eq:def-discrep1},
its finiteness can be proven analogously as in Theorem \ref{pro:stop-rule} (i).
Moreover, the following estimate holds true (compare with Item~(iii) of
Theorem \ref{pro:stop-rule})
\begin{equation} \label{eq:estim-log-L2}
d(y_j^\delta, \A_j x_{k_*^\delta}^\delta) \leq
\tau \delta_j \| \ln(y_{[k]}^\delta / (\A_{[k]} x_{k_*^\delta}^\delta) \|_{L^2} ,
\qquad j = 0, \dots, N-1 \,.
\end{equation}
\item In  Theorem \ref{pro:stop-rule} (iv), if one substitutes the assumption
\eqref{eq:assump-stabil-th} by $\| y_j - y_j^l \|_{L^2} \le \delta^l_j$,
$j = 0, \dots, N-1$, $l \in \mathbb N$, then the proof of the stability
result carries on with analogous argumentation.
\end{itemize}

Notice that the estimate in \eqref{eq:estim-log-L2} allows for the following
interpretation: The loping OS-EM iteration should be stopped at the index
$k_*^\delta$ (an integer multiple of $N$) when for the first time
\eqref{eq:estim-log-L2} is satisfied within a whole cycle.

The advantage of using this stopping rule resides on the fact that no
quantitative information on the constants $m$, $M$, $m_1$, $M_1$ is required
to compute the iteration. In other words, the constant $\gamma$ is
not required neither to test the ``loping condition" \eqref{eq:lop-cond-l2}
nor to verify the stopping rule based on \eqref{eq:estim-log-L2}.
This is obviously not the case if the ``loping condition"
$f_{[k]}^\delta(x_k^\delta) > \tau \gamma \delta_{[k]}$ in
\eqref{eq:def-osem-noise2} is to be implemented.
\end{remark}

\section{Numerical example} \label{sec:5}

In this section we compare the numerical performance of our loping
OS-EM method with the OS-EM and EM methods. As benchmark problem we
use a system of linear equations for the circular Radon transform.
The inversion of the circular Radon is relevant for the emerging
{\em photoacoustic computed tomography} \cite{KucKun08, PalNusHalBur07,
SchGraGroHalLen08, XuWan06}.

\medskip
Let $\epsilon < 1$ be some small positive number, let $\Omega :=
B_{1-\epsilon}(0) \subset \R^2$ denote the disc with radius $1-\epsilon$
centered at the origin, set
$$
    \Sigma_j :=  \left( \frac{2j\pi}{N}, \frac{2(j+1)\pi}{N} \right)
    \times (0,2) \, , \quad j = 0 , \dots,  N-1 \,,
$$
and let $ \Phi: \R \to \R$ be a continuous nonnegative function with
$\operatorname{supp} (\Phi) = [-\epsilon, \epsilon]$ and $\int_\R \Phi = 1$.

Our aim is the stable solution of \req{sys-ip}, with
$\A_j x  :=  \Phi \ast_r ( \Mo_j x)$, where
\begin{equation} \label{eq:smean}
   (\Mo_j x) (\ph, r)
   : =
   \frac{r N}{2\pi}
   \int_{S^1} x ( (\cos \ph, \sin \ph) + r \omega ) \, d \omega \,,
   \qquad  (\ph, r) \in \Sigma_j \,,
\end{equation}
is the \emph{circular Radon transform} restricted to $\Sigma_j$, and $\Phi \ast_r y = \Io_\Phi y$ denotes the
convolution of $\Phi$ and $y$ .
In \req{smean}, $x$ is considered as an element in $L^1(\R^2)$
by extending it with zero outside of $\Omega$.

One  verifies that the operators $\A_j$ can be written in the form
(\ref{eq:sys-ops}), with $s =( \ph, r)$ and
$$
    a_j( t, \ph, r) =
    \Phi
    \left(
    \abs{( \cos \ph, \sin \ph ) - t} - r
    \right)\,, \quad  j = 0, \dots, N-1 \,.
$$
Moreover, the adjoint of $\A_j$ is given by $\A_j^* y =
\Bo_j  (\Phi \ast_r y)$, where
\begin{equation}\label{eq:bp}
    ( \Bo_j y)(t)
    =
    \frac{N}{2\pi}
    \int_{2j\pi/N}^{2(j+1)\pi/N}
    y \bigl( \abs{t-(\cos \ph, \sin \ph)} \bigr) d \ph \,,
\end{equation}
is the \emph{circular backprojection}.
Hence $\A_j^* 1 = 1$ and the operators $\A_j$ satisfy assumption (A1).
However, since $a_j$ are not bounded from below, $\A_j$ do not satisfy (A2).

\begin{remark}
For any positive $\lambda$, the operators
$$
    \A_j^{(\lambda)} x := \frac{1}{1 + \lambda\abs{\Sigma_j}}
                         \left( \A_j x+ \lambda \int_\Omega x \right)
\,, \qquad j=0, \dots, N-1\,,
$$
clearly satisfy (A2). Since $(\A_j^{(\lambda)})^* y =
( \A_j^* y + \lambda \int_{\Sigma_{j}} y ) / (1 + \lambda\abs{\Sigma_j})$ we have $(\A_j^{(\lambda)})^* 1 = 1$, proving that (A1) is also satisfied.

Therefore, we shall consider for the rest of this section the
system of equations
\begin{equation} \label{eq:iter-mod}
    \A_j^{(\lambda)} x = y_j^{(\lambda)}
    :=
    \frac{1}{1 + \lambda\abs{\Sigma_j}}\left( y_j + \lambda \int_{\Sigma_{j}} y_j \right)\, ,
    \qquad j = 0, \dots, N-1 \,.
\end{equation}
The identity  $\int_\Omega x = \int_\Omega x \A_j^* 1 = \int_{\Sigma_{j}} \A_j x$ implies that
$x$ is a solution of \eqref{eq:iter-mod} if and only if $x$
satisfies  $\A_j x = y_j$.
\end{remark}

If noisy data $y_j^\delta$ with
$\norm{y_j^{\delta} - y_j}_{L^1} \leq \delta_j$ are
available, then
\begin{equation*}
\norm{y_j^{(\lambda), \delta} - y_j^{(\lambda)}}_{L^1}
\leq
\frac{1}{1 + \lambda\abs{\Sigma_j}}
\left(
\norm{y_j^{\delta} - y_j}_{L^1} +
\lambda \Bigl| \int_{\Sigma_{j}} y_j^{\delta}-y_j \Bigr|
\right)
\leq
\frac{\delta_j(1+\lambda)}{1 + \lambda\abs{\Sigma_j}}\,,
\end{equation*}
where $y_j^{(\lambda), \delta}$ is defined in the same way as
$y_j^{(\lambda)}$, with $y_j$ replaced by $y^\delta_j$.
Therefore, the loping OS-EM iteration with noisy data $y^\delta_j$
applied to system \eqref{eq:iter-mod} reads as
\begin{equation} \label{eq:losem-la}
\begin{aligned}
x_{k+1}^\delta
& :=
x_k^\delta \, \omega_k \,,
\\
\omega_k
& :=
  \begin{cases}
   \frac{\Bo_{[k]} \Io_\Phi  + \lambda}{1+\lambda\abs{\Sigma_j}}
   \Bigl(
   \frac{y_{[k]}^{(\lambda), \delta}}
   {\Io_\Phi\Mo_{[k]} x_k^\delta + \lambda}
   \Bigr)
   \,,
   &  \frac{d( y_{[k]}^{\delta} + \lambda,
   \Io_\Phi \Mo_{[k]} x_k^\delta + \lambda)}{1 + \lambda} > \tau \gamma
   \delta_{[k]}\,, \\
   1\,, &  \text{ else}\,.
  \end{cases}
  \end{aligned}
\end{equation}
Here we made use of the fact that the initial guess satisfies $\int_\Omega x_0^\delta = 1$,
which implies $\int_\Omega x_k^\delta = \int_{\Sigma_{[k]}} \A_{[k]} x_k^\delta = 1$
for every $k$.

\begin{remark}
Iteration \req{losem-la} assumes continuous data $y^{\delta}_j \in L^1(\Sigma_j)$, whereas
in practical applications only discrete data are available.
In the following we assume that data
\begin{equation*}
    \yd_j^{\delta}[ i_\ph, i_r ]
    :=
    y_j^{\delta} ( \phd[i_\ph],  \rd[i_r]  )
    \,,
    \qquad(i_\ph, i_r) \in   \set{ j N_\ph, \dots,  (j+1) N_\ph -1} \times \set{ 0, \dots,  N_r}\,,
\end{equation*}
are given, with $\phd[i_\ph] := 2 i_\ph \pi / N_\ph$, $\rd[i_r] := 2 i_r /N_r$, and
$N_\ph$, $N_r+1$ denoting the number of samples of $y_j^\delta$
in the angular and radial variable, respectively.
\end{remark}

In the numerical implementation $\Mo_j$, $\Bo_j$, $\Io_\Phi$
and $d$ are replaced (as described below) with finite dimensional approximations $\Md_j$, $\Bd_j$, $\Id_\Phi$, $\f d$,
and \eqref{eq:losem-la} is approximated by
\begin{equation}\label{eq:losem-d}
\begin{aligned}
x_{k+1}^\delta(\td[i])
\simeq
\xd_{k+1}^\delta[i]
&:=
\xd_k^\delta[i] \omd_k[i] \,, \qquad i \in \set{0, \dots, N_t}^2
\\
\omd_k
&:=
 \begin{cases}
 \frac{\Bd_{[k]}\Id_\Phi + \lambda}{1+4\pi\lambda/N}
   \Bigl(
   \frac{\yd_{[k]}^{\delta} +\lambda }
   {\Id_\Phi \Md_{[k]} \xd_k^\delta + \lambda}
   \Bigr)
     \,,
   &  \frac{\f d\bigl(\yd_{[k]}^{\delta}+ \lambda, \Id_\Phi \Md_{[k]}^{(\lambda)} \xd_k^\delta + \lambda \bigr)}{1 + \lambda} > \tau \gamma \delta^{(\lambda)}_{[k]}\,,
   \\
   1\,, &  \text{ else}\,.
  \end{cases}
  \end{aligned}
\end{equation}
Here $\yd_j^{\delta} = (\yd^{\delta}_j[ i_\ph, i_r ])_{i_\ph, i_r}$,
$\xd_k^\delta = (\xd_k^\delta [i])_i$, $\omd_k = (\omd_k[ i])_i$,
and $\td[i] = -(1,1) + 2i/N_t$ with $(N_t+1)^2$ denoting the
number of samples in the variable $t$.

\begin{enumerate}
  \item
  The discretized circular Radon transform
  \[\Md_j : \R^{(N_t+1) \times (N_t+1)} \to \R^{N_\ph \times (N_r+1)}\]
  is obtained by replacing $x$ in \req{smean} with the bilinear spline $T(\xd)$ satisfying
  $T(\xd)(\td[i]) = \xd[i]$, and approximating the resulting integrals over the $S^1$
  with the trapezoidal rule. This leads to
  \begin{equation}\label{eq:circ-d}
  (\Md_j \xd) [i_\ph, i_r]
  =
  \frac{2}{N_t}
  \sum_{i_\omega = 0}^{3 \rd[i_r] N_t}
  T(\xd)
  \Bigl(
  \sid[i_\ph] + \rd[i_r] \omd[i_\omega]
  \Bigr) \,,
  \end{equation}
  where $\sid[i_\ph] := (\cos \phd[i_\ph], \sin\phd[i_\ph])$, $\omd[i_\omega] := \bigl( \cos(2\pi i_\omega/N_t ),\sin(2\pi i_\omega/N_t) \bigr)$, and $3 \rd[i_r] N_t$ is the number of supporting
  points when applying the trapezoidal rule.

  \item
  Assuming that $\epsilon = 2 K /N_r$
  for some $K \in \mathbb N$,
  the convolution $\Io_\Phi y  = \Phi \ast_r y$
  is approximated by
  \[
  (\Id_\Phi \yd)[i_\ph, i_r]
  =
  \frac{2}{N_r}
  \sum_{i_r' = i_r-K}^{i_r+K}
  \Phi \bigl( 2(i_r'-i_r)/N_r \bigr)
  \yd[i_\ph, i_r'] \,,
  \]
  where $\yd[i_\ph, i_r] := 0$
  for $i_r$ outside $\set{0, \dots, N_r}$.

\item
The discretized back-projection
$\Bd_j: \R^{N_\ph \times (N_r+1)} \to \R^{(N_t+1) \times (N_t+1)}$ is  defined by
\begin{equation*}
    (\Bd_j \yd)[i]
    :=
    \frac{N}{N_\ph}
    \sum_{i_\ph = j N_\ph}^{(j+1)N_\ph-1}
        T_r(\yd)
        \bigl(i_\ph, \abs{ \td[i] - \sid[i_\ph] }
        \bigr) \,,
\end{equation*}
if $\td[i]  \in \Omega$, and setting
$(\Bd_j \yd)[i] := 0$ for $\td[i] \not \in \Omega$.
Here $T_r(\yd) $ denotes  the piecewise linear spline in the second variable satisfying
$T_r(\yd) \bigl( i_\ph,  \rd[i_r] \bigr) = \yd[i_\ph, i_r]$.

\item
Finally, the discrete KL-distance is defined by
\[
\f d( \f v, \f u)
=
\frac{4\pi}{N_r N_\ph}\sum_{i_\ph = jN_\ph}^{(j+1) N_\ph - 1}
\sum_{i_r = 0 }^{N_r}
\f v[i_\ph, i_r]
\ln \frac{\f v[i_\ph, i_r]}{ \f u[i_\ph, i_r]}
-
\f v[i_\ph, i_r]
+
\f u [i_\ph, i_r] \,.
\]
for $\f v, \f u \in \R^{N_\ph \times (N_r+1)}$.
\end{enumerate}

\begin{remark}[Numerical Complexity]
Assuming $N_t = N_r$, the numerical complexity for performing one iteration cycle
(which consists of $N$ subsequent steps in \req{losem-d}) is $\mathcal O (N_{\rm angle}  N_t ^2)$.
Here  $N_{\rm angle} = N_\ph N$ corresponds to the overall
angular data samples, which is independent of $N$ in practice.
Therefore in the following we always compare the reconstruction error in dependence of the number of iteration cycles.
\end{remark}

\begin{figure}[htb!]
  \includegraphics[width = 0.48\textwidth]{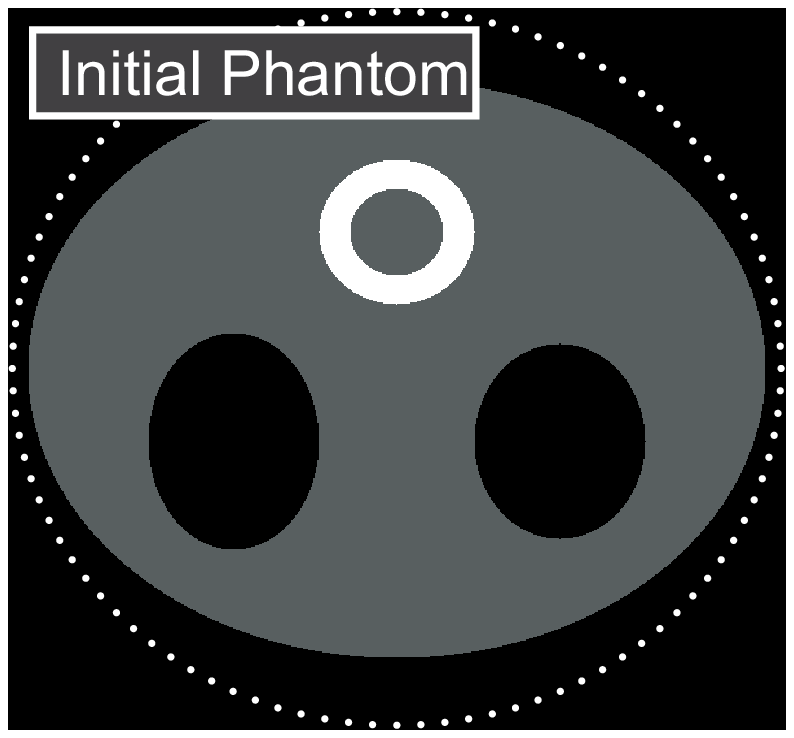}\quad
  \includegraphics[width = 0.48\textwidth]{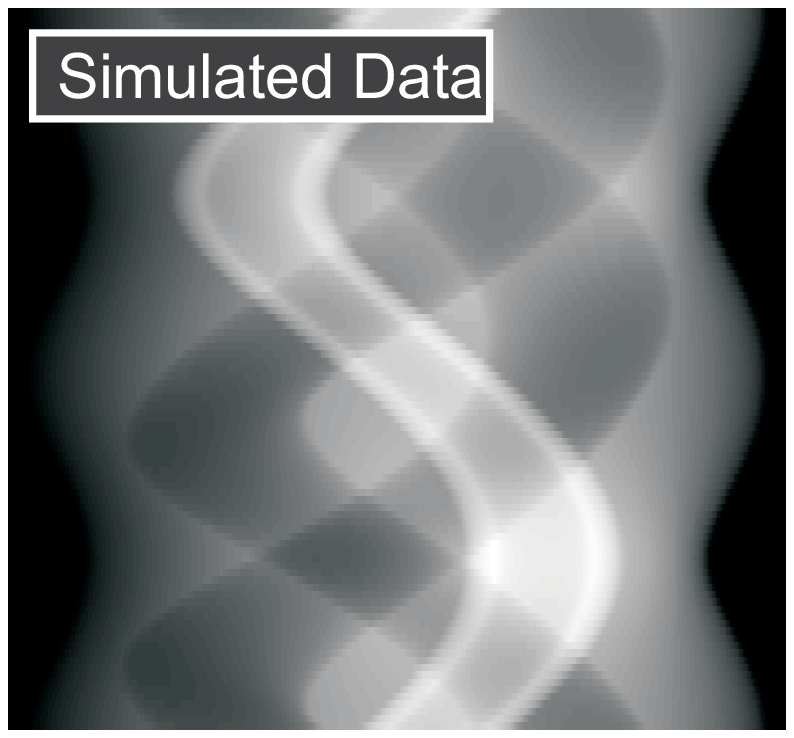}
  \caption{Original phantom $x^*$ (left) and simulated  data $(\A_j x^*)_j$.}\label{fig:phant-data}
\end{figure}

In the following numerical examples we apply the (loping) OS-EM
iteration with $N = 1$ (corresponding to the EM algorithm), $N = 5$,
$N=10$, and $N=20$ subsets. The original phantom $x^*$ (the exact
nonnegative solution) is shown in the left picture in
Figure~\ref{fig:phant-data} and consists of a superposition of
characteristic functions. Note that a similar phantom was
reconstructed in \cite{HudLar94} where the OS-EM technique was
introduced. The data $\yd_j$, shown in the right picture in
Figure~\ref{fig:phant-data}, were calculated numerically by \req{circ-d}
for $N_{\rm angle} = N_\ph N =100$ angular samples. In order to avoid inverse
crimes, much larger $N_t$ is used for the data simulation as for the
application of the loping OS-EM iteration. In all examples $\xd_0 =
\xd_0^\delta = 1/\bigl((1-\epsilon)^2\pi\bigr)$ is used as initial
guess and the parameters $\epsilon$ and $\lambda$ are chosen to be
$0.02$ and $0.01$, respectively.

\begin{figure}[htb!]
  \includegraphics[width=0.32\textwidth]{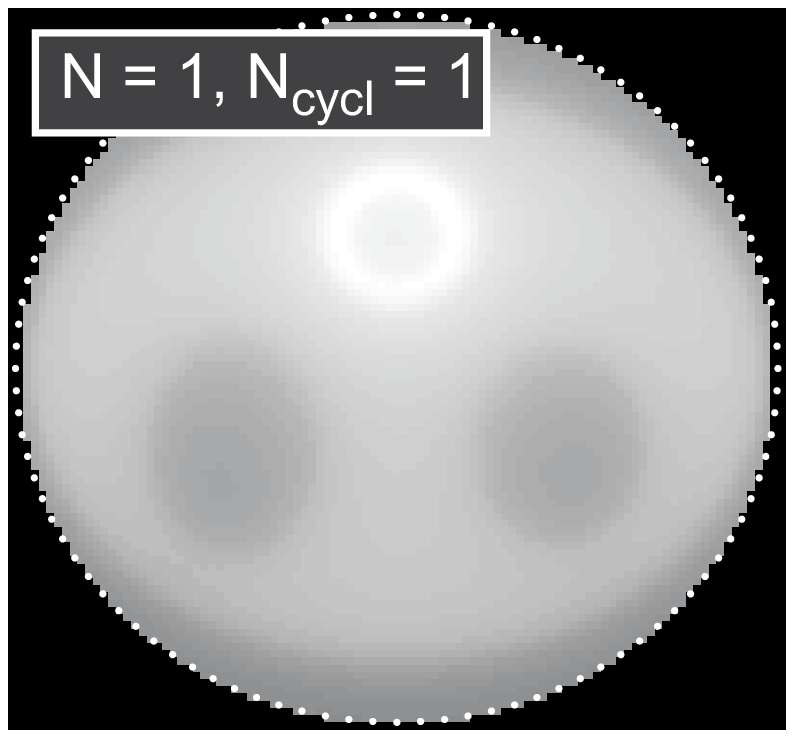}
  \includegraphics[width=0.32\textwidth]{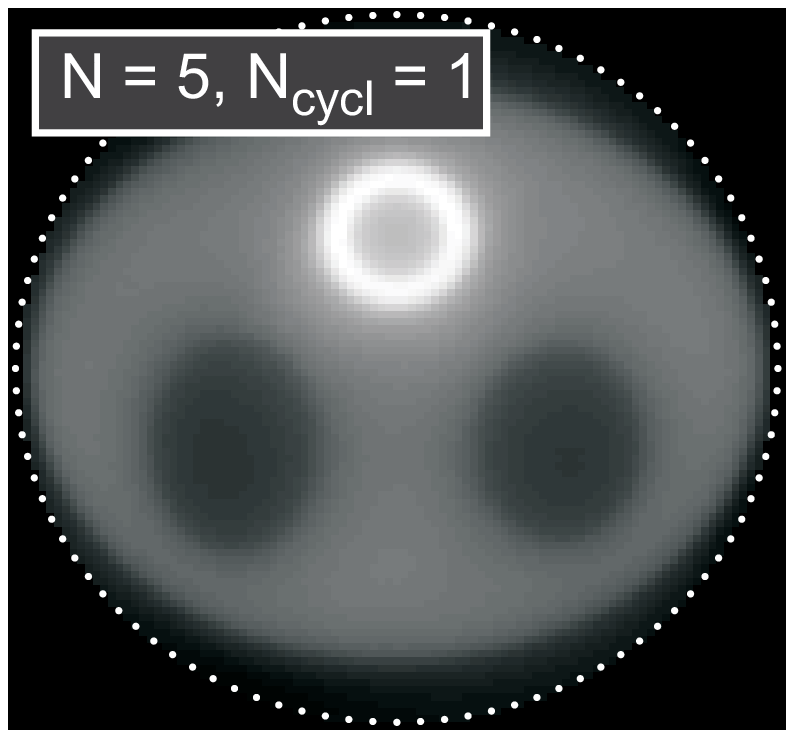}
  \includegraphics[width=0.32\textwidth]{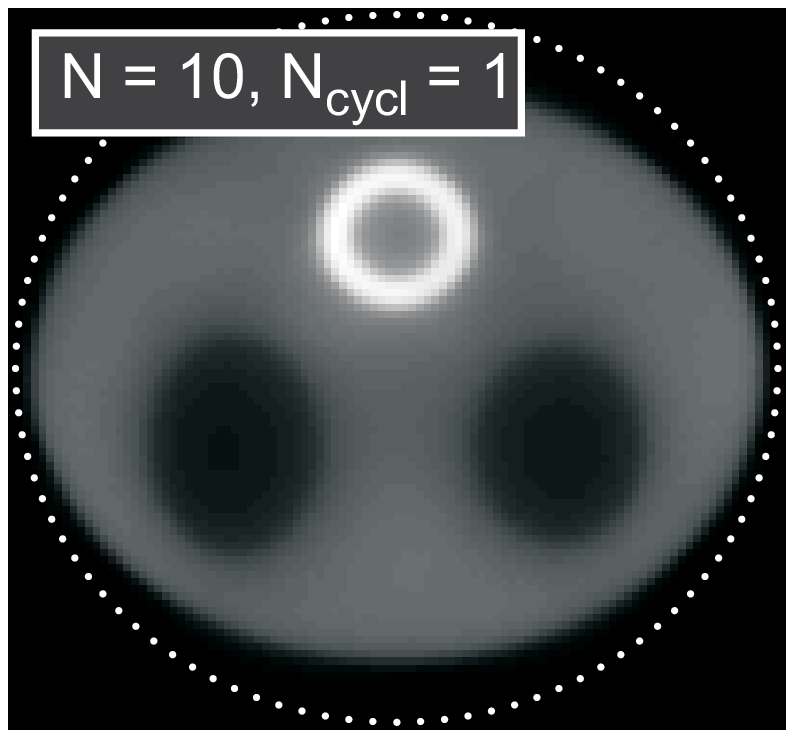}\\[0.4ex]
  \includegraphics[width=0.32\textwidth]{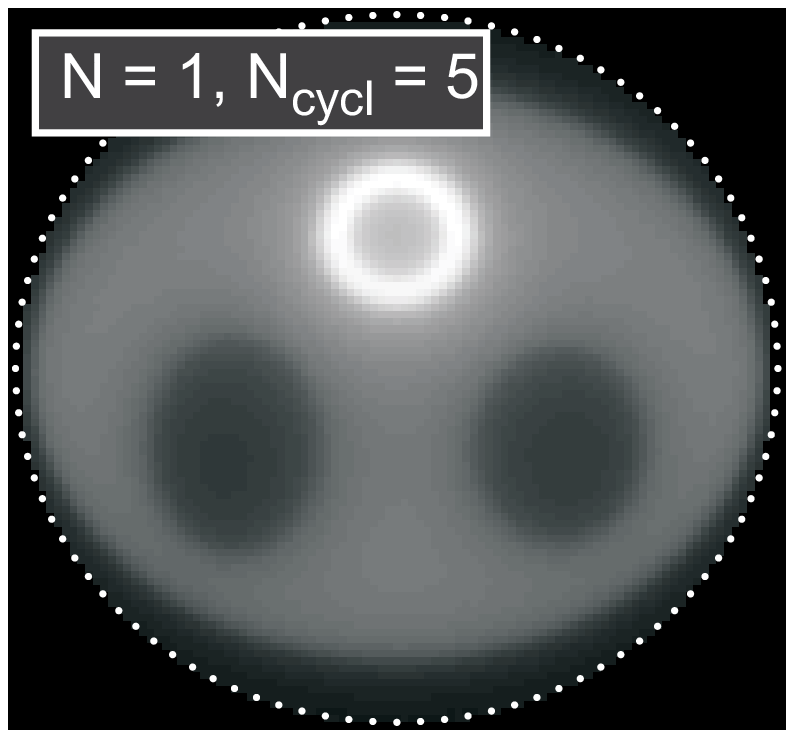}
  \includegraphics[width=0.32\textwidth]{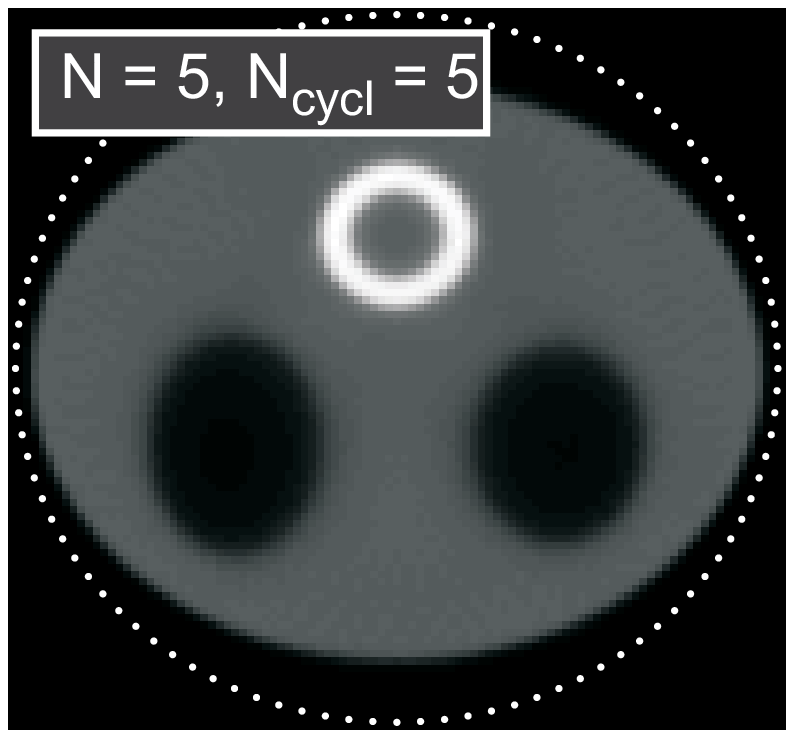}
  \includegraphics[width=0.32\textwidth]{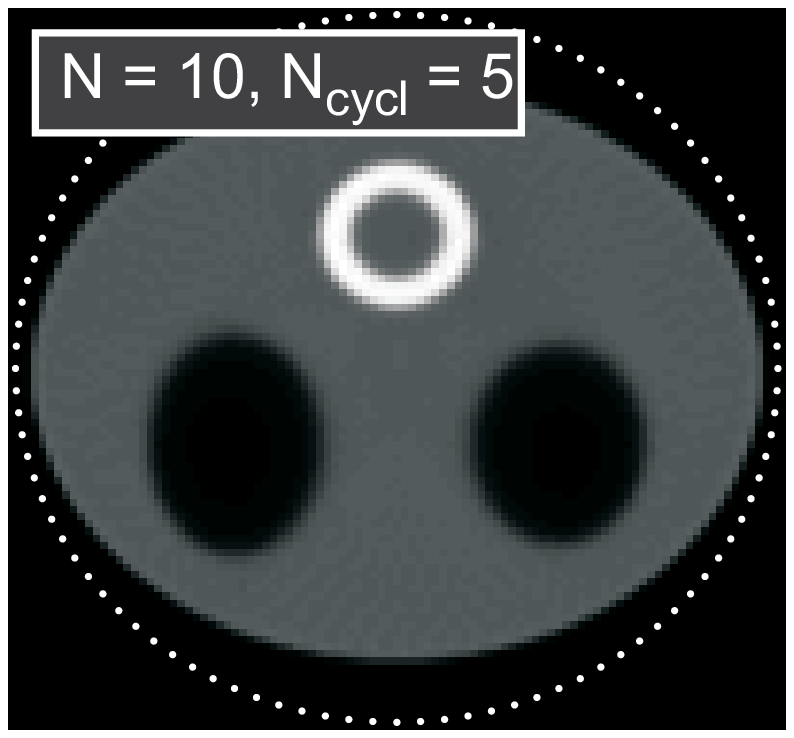}\\[0.4ex]
  \includegraphics[width=0.32\textwidth]{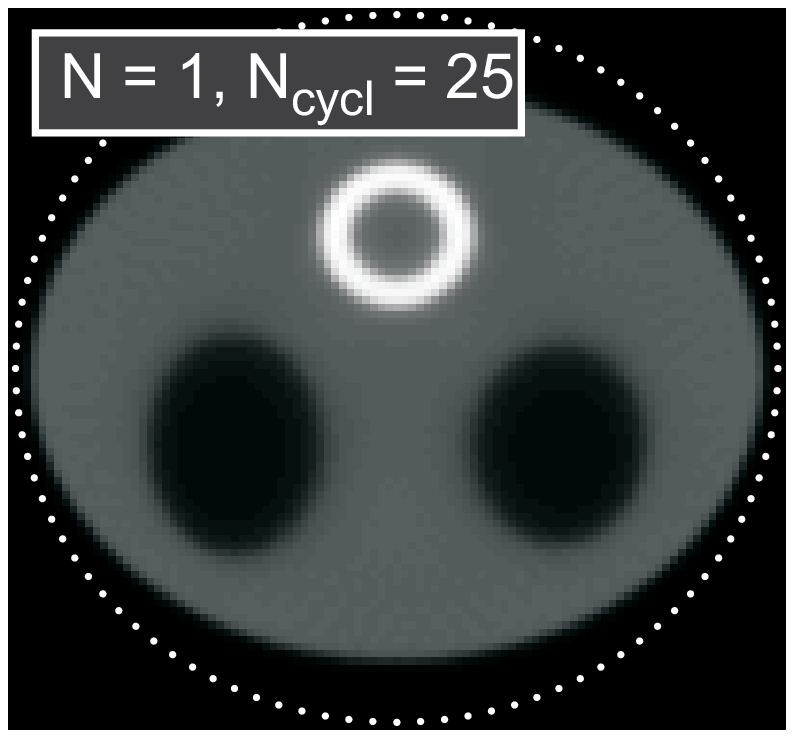}
  \includegraphics[width=0.32\textwidth]{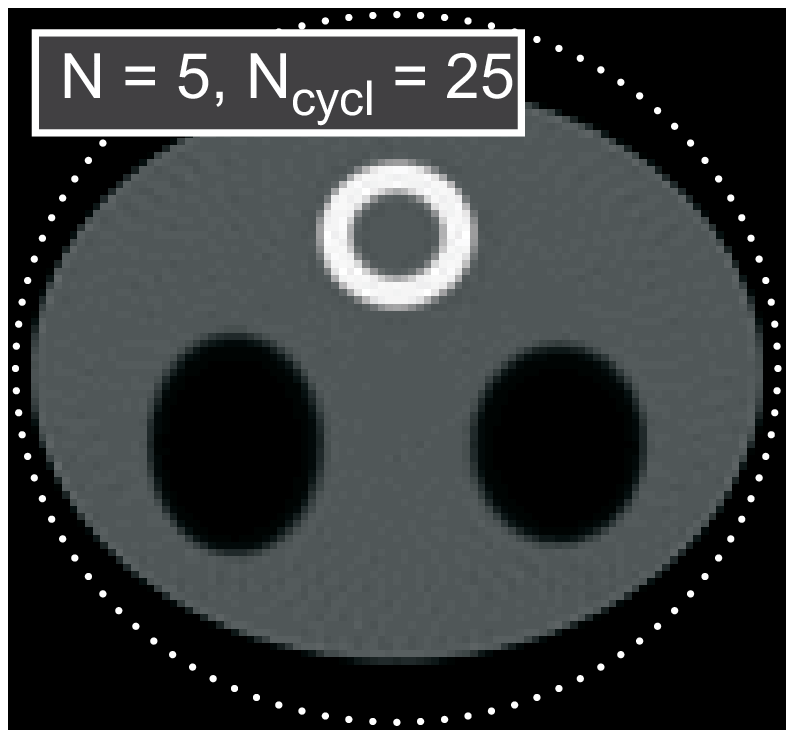}
  \includegraphics[width=0.32\textwidth]{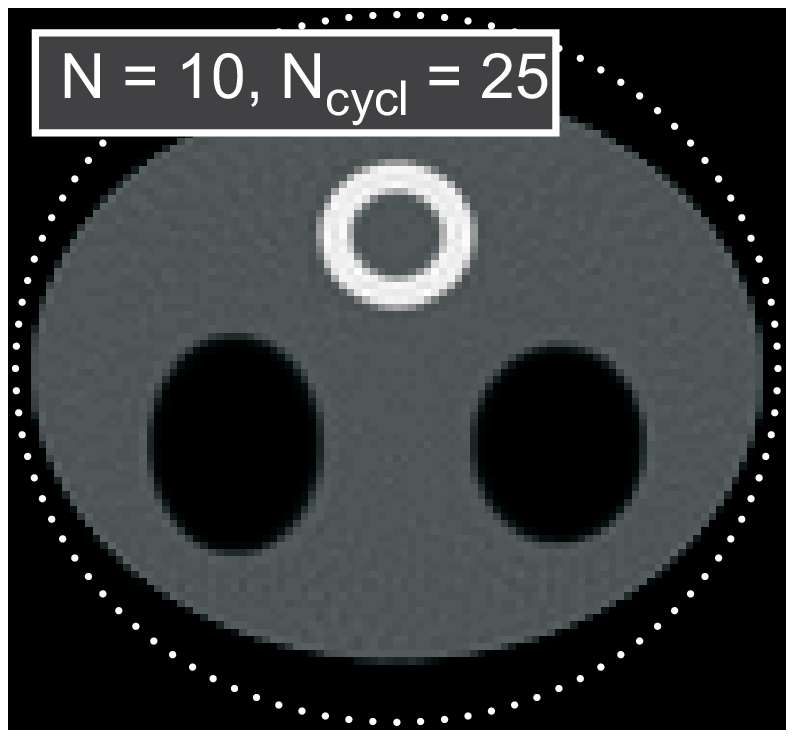}
  \caption[Iterates for exact data]{Exact data experiment: Iterates for $N_t = N_r=100$ with $N=1$ (left), $N=5$ (middle)
  and $N=10$ (right) after $1$, $5$ and $25$ cycles.}
  \label{fig:exact}
\end{figure}

\begin{psfrags}
\psfrag{Klerr100}{$\log \f d( \xd^*, \xd_k)$ for $N_t=100$}
\psfrag{Klerr200}{$\log \f d( \xd^*, \xd_k)$ for $N_t=200$}
\begin{figure}[htb!]
  \includegraphics[width=0.44\textwidth]{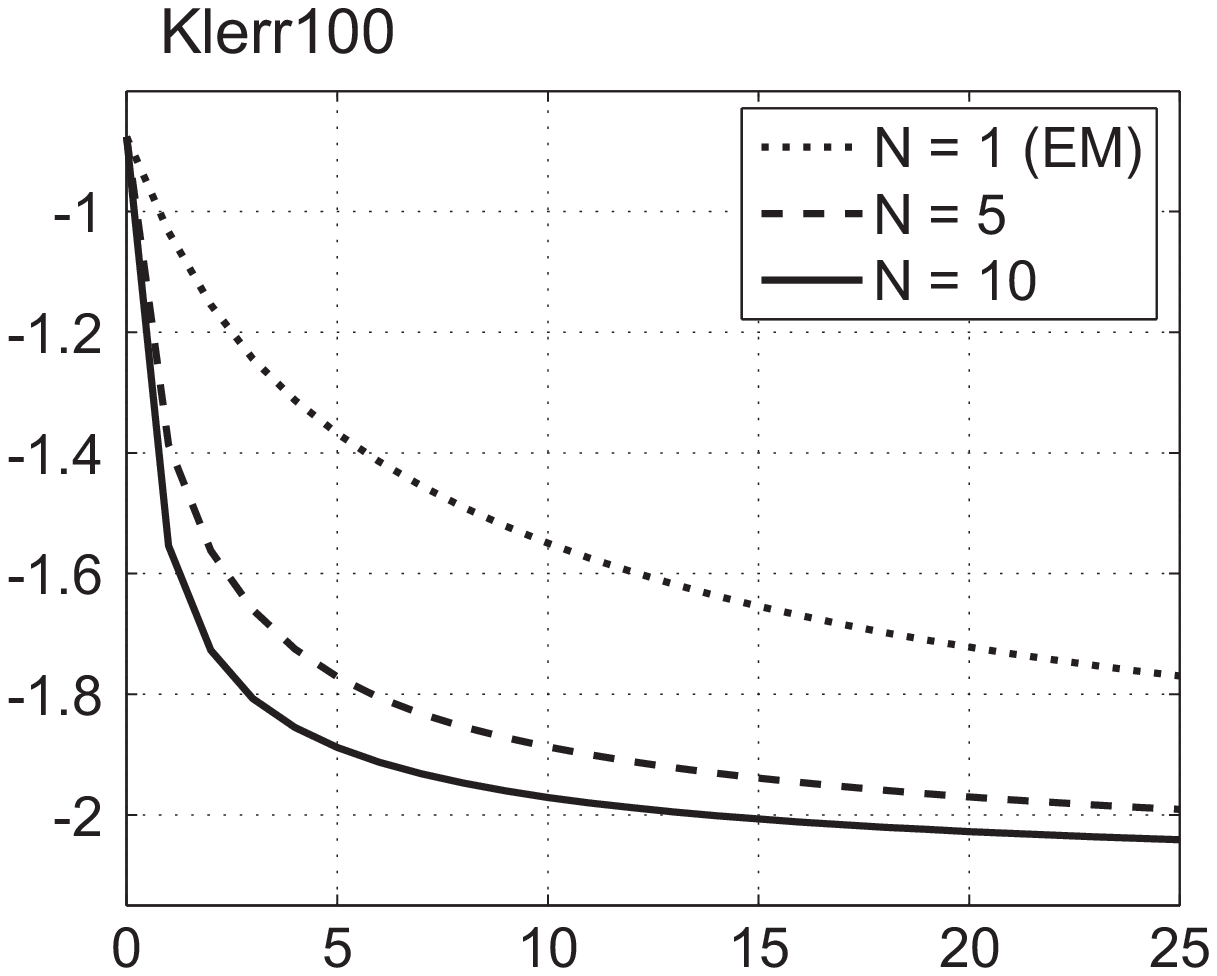}\qquad
  \includegraphics[width=0.44\textwidth]{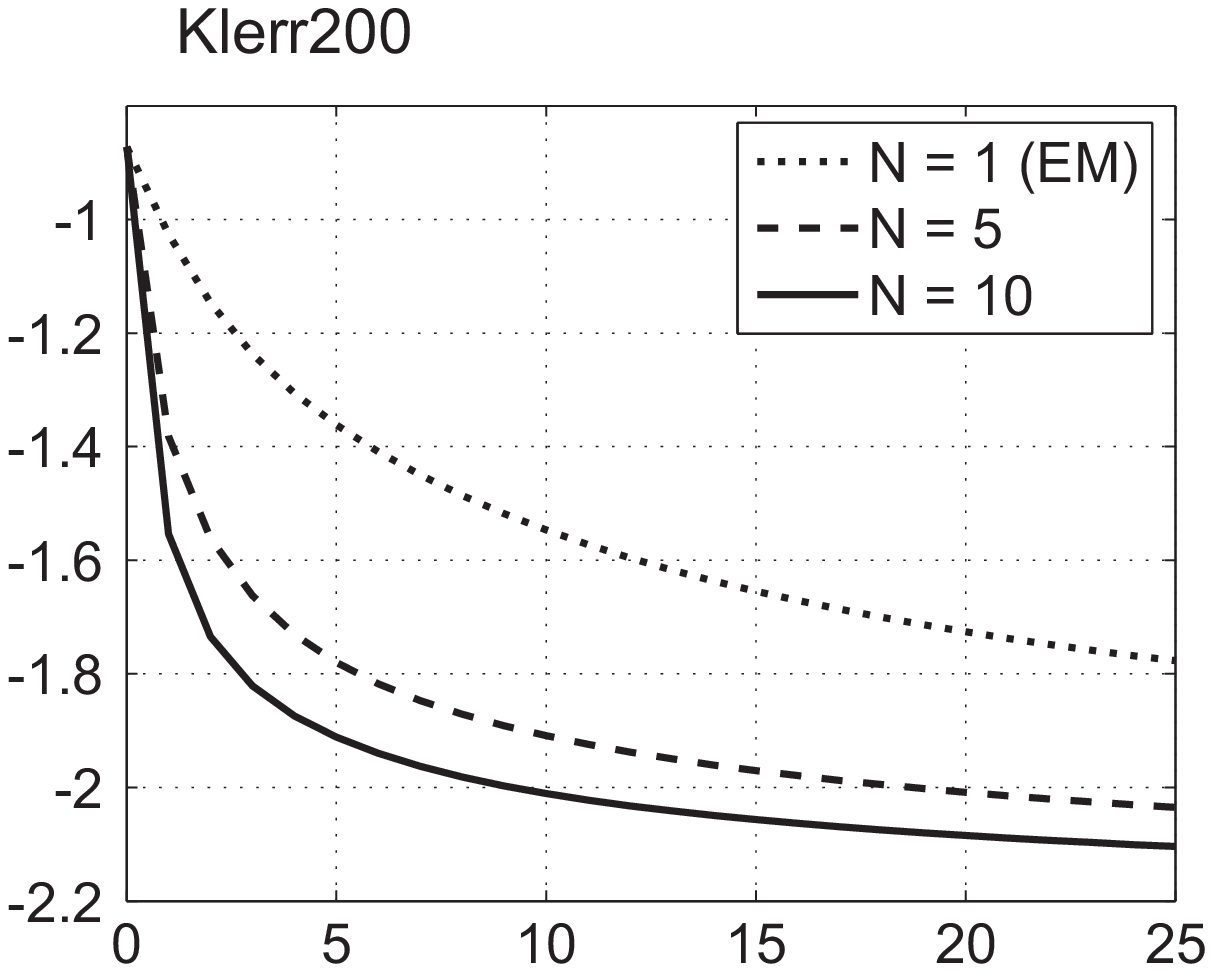}
  \caption[Error plots for exact data]{Exact data experiment: Logarithmic plots of iteration error with
  respect to the Kullback-Leibler distance for $N_t=N_r=100$ (left) and $N_t=N_r=200$ (right).}
  \label{fig:exact-error}
\end{figure}
\end{psfrags}

The iterations of the OS-EM method applied to exact data with $N_t =
N_r = 100$ and different values of $N$ are depicted in
Figure~\ref{fig:exact}. It can be seen that the 5-th iterate with EM
has similar quality as the 1-th iterate  with OS-EM for $N = 5$. As
a rough rule one can say that making $N$ cycles with the EM
algorithm leads to an improvement similar to $1$ cycle with the
OS-EM algorithm. This can also be recognized in
the left image in Figure~\ref{fig:exact-error}, where the evolution of the error is
depicted with respect to the KL-distance.

In order to investigate the dependence of the OS-EM iteration on the discretization level,  we repeated the experiment with
$N_t = N_r = 200$.
The right image in Figure~\ref{fig:exact-error} shows
the corresponding logarithmic error. As expected, the error is relatively independent
on the discretization.

\medskip
In the case of noisy data we apply the loping OS-EM iteration
\req{losem-d}. The noisy data $\yd_j^\delta$ is created by adding
$5\%$ Poisson distributed noise to the simulated data $\yd_j$, such that
$4\pi/(N_r N_\ph) \sum \bigl|\yd_j[i_\ph, i_r] - \yd_j^\delta[i_\ph, i_r]\bigr|
\cong 0.05$.

\begin{remark}
Our numerical experiments show that, for large $\delta$ and $ \tau \cong 1$,
far too many iterations are loped.
A significant improvement can be obtained if $\tau = \tau(\delta)$ is chosen
in dependence of the noise level, with $\tau(\delta) < 1$ for large $\delta$
and $\tau(\delta)$ converging to some $\tau_\infty > 1$ for $\delta \to 0$.
It is clear that the asymptotic convergence analysis (for $\delta \to 0$)
remains valid in such a situation.
\end{remark}

\begin{figure}[htb!]
\centering
  \includegraphics[width=0.45\textwidth]{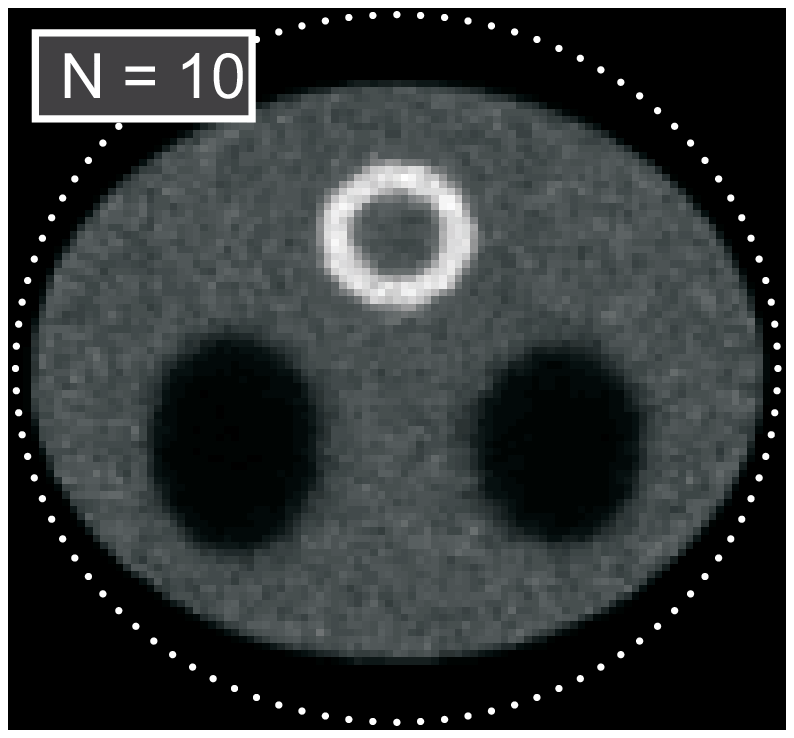}
  \includegraphics[width=0.45\textwidth]{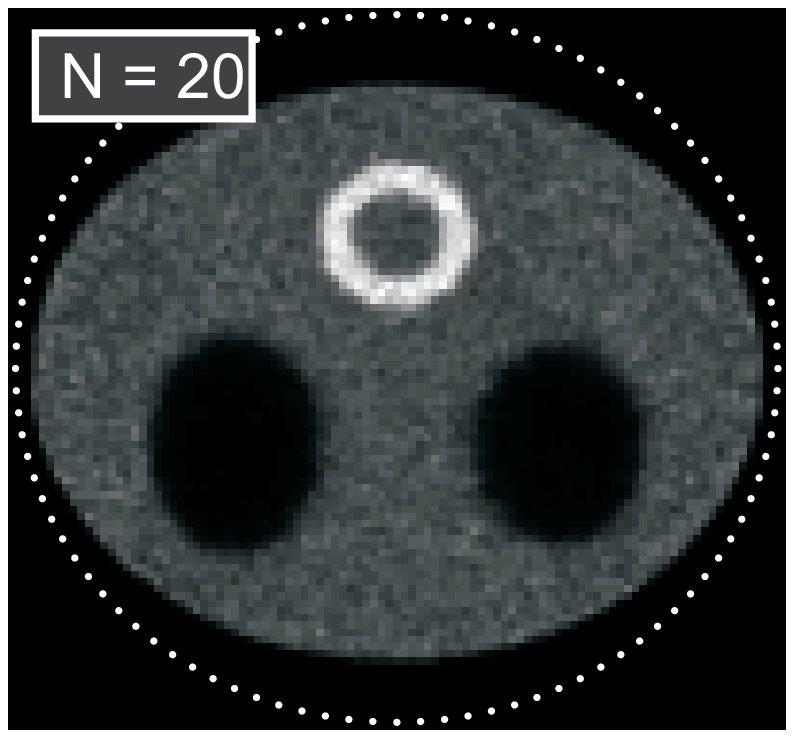} \\[0.4ex]
  \includegraphics[width=0.45\textwidth]{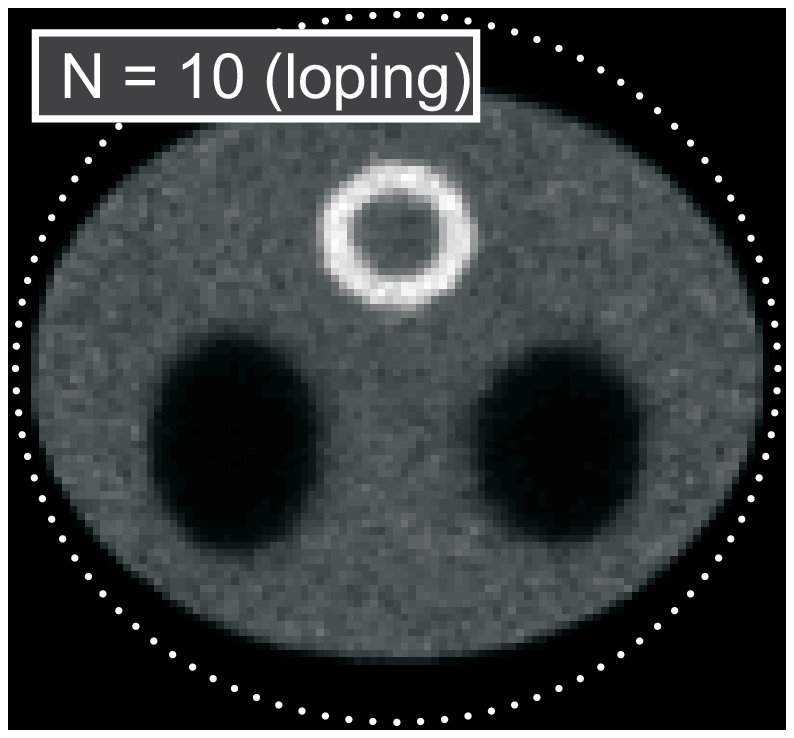}
  \includegraphics[width=0.45\textwidth]{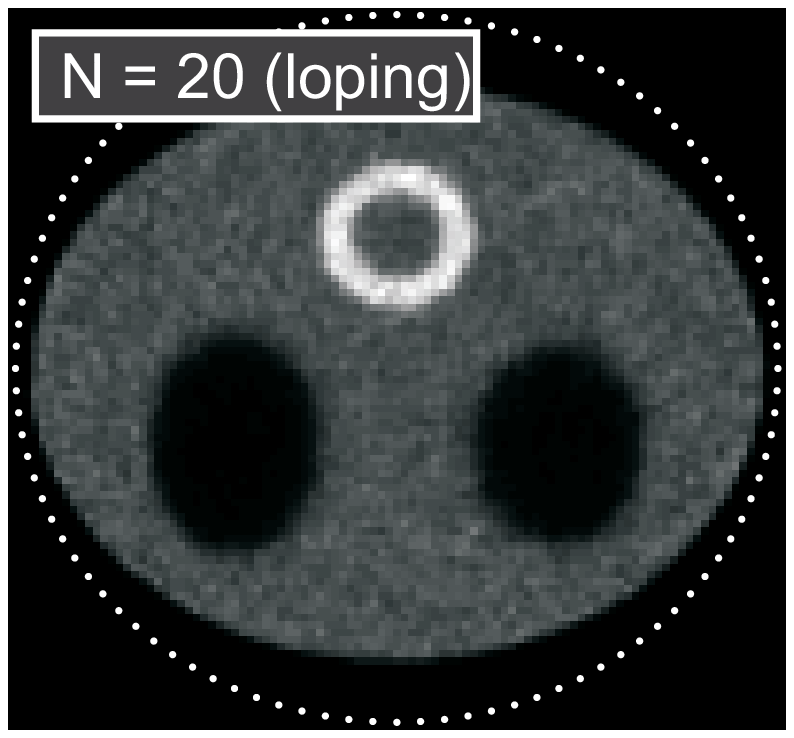}
  \caption[Reconstructions from noisy data]{Noisy data experiment: Iterates without loping (top line) and
  with loping (bottom line). The loping iterations are stopped
  automatically whereas their non-loping counterparts are stopped at the
  iteration cycle where $d(x^*, x_k^\delta)$ is minimal (which is not
  available in practice).} \label{fig:noisy}
\end{figure}

\begin{psfrags} \begin{figure}[htb!]
  \begin{center}
  \includegraphics[width = 0.5\textwidth]{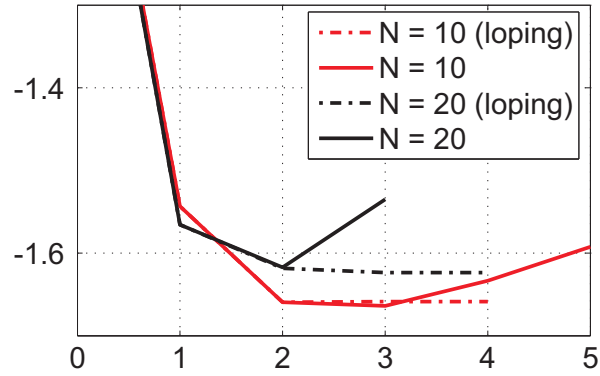}
  \caption[Error plots for noisy data]{Noisy data experiment:
  Evolution of the relative error $\log \f d( \xd^*, \xd_k^\delta)$
  for loping and non-loping OS-EM iterations.} \label{fig:noisy-error}
  \end{center}
\end{figure}

\begin{figure}[htb!]
  \begin{center}
  \includegraphics[width = 0.45\textwidth]{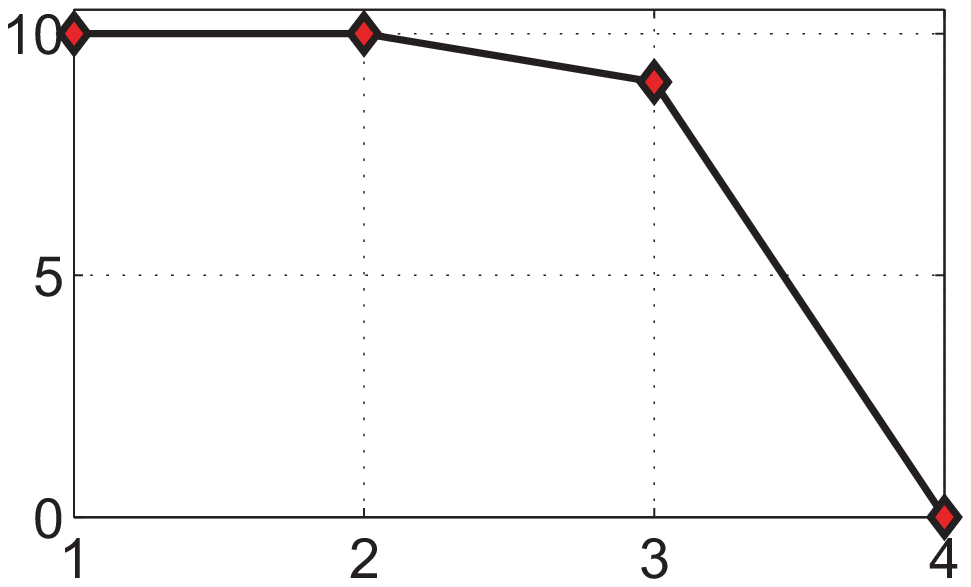}\qquad
  \includegraphics[width = 0.45\textwidth]{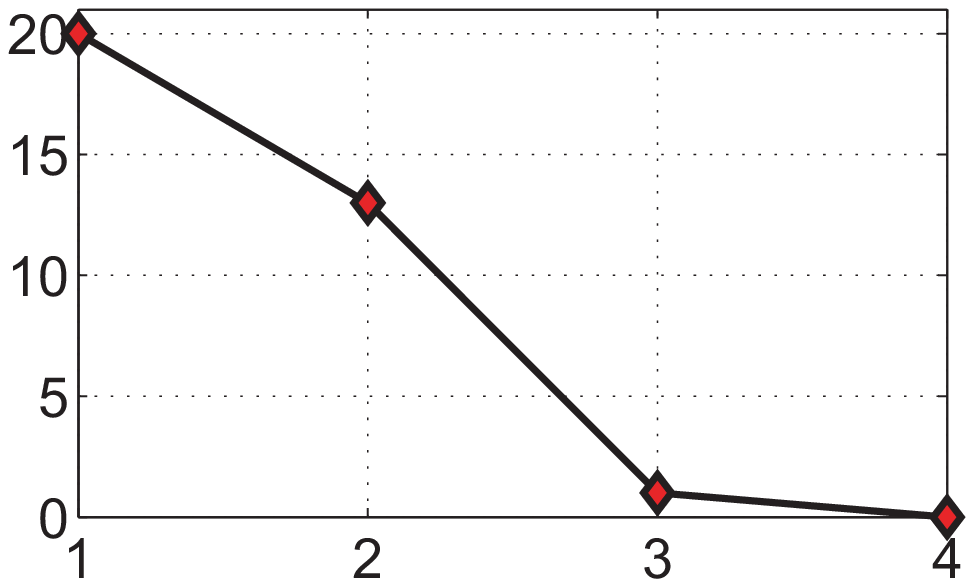}
  \caption[Number of performed cycles]{Noisy data experiment: The $x$-axis shows the number of cycles,
  while the number of actually performed iterations within each cycle is
  shown at the $y$-axis.} \label{fig:steps}
  \end{center}
\end{figure} \end{psfrags}

The reconstruction  for noisy data with $N_t=N_r=100$ are depicted
in Figure~\ref{fig:noisy}. For comparison purposes, results of the
OS-EM iteration (without loping strategy) are also included. The
loping OS-EM is automatically stopped according to
\eqref{eq:def-discrep1} whereas their non-loping counterparts are
stopped after the cycle with minimal error $\f d(\xd^*, \xd_k^\delta)$,
which is not available in practice. All reconstructions are quite
comparable. Figure~\ref{fig:noisy-error} shows the evolution of the
error with respect to the KL-distance. In this figure one also
notices the semi-convergence of the non-loping iterations, which
happens typically when applying non-regularized iterative schemes to
ill-posed problems \cite{DeuEngSch98, HanNeuSch95, KalNeuSch08}.

\begin{table}[tbh!]
\centering
\begin{tabular}{l | @{\hspace{1cm}} c @{\hspace{0.5cm}} c @{\hspace{0.5cm}}
                c @{\hspace{0.5cm}} c @{\hspace{0.5cm}} c }
\toprule
& $N$ & $N_{\rm cycl}$ & time (sec) & $\f d(\xd^*, \xd_k^\delta)$ \\
\midrule
loping OS-EM   & 10    &  4   &  13.4   & 0.022 \\
OS-EM       & 10    &  3   &  9.2    & 0.022 \\
loping OS-EM   & 20    &  4   &  13.4   & 0.024 \\
OS-EM       & 20    &  2   &  6.3    & 0.024 \\
\bottomrule
\end{tabular}
\caption[Numerical performance]{Comparison of the performance of different iterative methods.
The non-loping iterations are stopped after the cycle with minimal error,
whereas the loping OS-EM are automatically stopped according to
\eqref{eq:def-discrep1}.} \label{tab:summary}
\end{table}

An inspection of Figure~\ref{fig:noisy-error} shows that the
regularized solution of the loping OS-EM methods (automatically
stopped) have errors comparable to the optimal solution of their
non-loping counterparts when stopped after the cycle with minimal
error (which is not available in the practice).
Figure~\ref{fig:steps} shows the number of actually performed
iterations. Table~\ref{tab:summary} summarizes run times and errors
with $N_t=N_r=100$, $N_{\rm angle} = 100$ (with non-optimized Matlab
implementation on HP Notebook with 2 GHz Intel Core Duo processor).

\section{Conclusions} \label{sec:6}

This article is devoted to the investigation of OS-EM type algorithms
for solving systems of linear ill-posed equations.
We focus on showing regularization properties of the
proposed methods.

In the case of exact data, our approach originates an algorithm analog to the OS-EM iteration.
We are able to prove  monotonicity results with respect
to the Kullback-Leibler distance as well as weak convergence in case of boundedness of the iterations.
In the noisy data case, we propose a loping OS-EM iteration which differs from the OS-EM method due to the introduction of a loping strategy.
This loping strategy renders the proposed iteration a regularization method.
We prove monotonicity of the iterates and study  stability properties of our method.

What concerns numerical effort, we conjecture that the loping OS-EM
algorithm is at least as efficient
as the well established OS-EM
method. The numerical experiments with \req{losem-d}  for inverting
the circular Radon transform support this conjecture. In the case of
exact data,  \req{losem-d} reduces to a discretized version of the
continuous OS-EM iteration applied to the system \req{iter-mod}.
However it is slightly different to the discrete OS-EM iteration of
\cite{HudLar94} since $\Bd_j$ is not the exact transpose of $\Md_j$.
Moreover, opposed to \cite{HudLar94} our continuous convergence
analysis applies independent on the discretization level.

\section*{Acknowledgments}
M.H. has been supported by the technology transfer office of the
University Innsbruck (transIT) within the framework of the NFN
``Photoacoustic Imaging in Biology and Medicine'' from the Austrian
Science Foundation (project S10505-N20). The work of A.L. is
supported by the Brazilian National Research Council CNPq, grants
306020/2006-8 and 474593/2007-0. E.R. acknowledges support from the
Austrian Science Foundation, Elise Richter scholarship (V82-N18
FWF).

\makeatletter
\@ifundefined{href}{\newcommand{\href}[2]{#2}\@latex@warning@no@line{href
  undefined, please load hyperref}}{}
\makeatother
\def\cprime{$'$} \def\cprime{$'$} \def\cprime{$'$}

%

\end{document}